\documentclass[prd,referee, nofootinbib]{revtex4}
\usepackage{mathrsfs}
\usepackage{yfonts}
\usepackage{tipa}
\usepackage{bm}
\usepackage{bbm}
\usepackage{amsmath}
\usepackage{amssymb}
\usepackage{amsthm}
\usepackage{amsfonts}
\usepackage{mathtools}
\usepackage{latexsym}
\usepackage{relsize}
\usepackage[greek,english]{babel}
\usepackage{booktabs}
\usepackage[iso-8859-7]{inputenc}
\usepackage{tikz}
\usepackage{array}

\def\be{\begin{equation}}
\def\ee{\end{equation}}
\def\bea{\begin{eqnarray}}
\def\eea{\end{eqnarray}}

\begin{document}

\title{Solving the nonlinear biharmonic equation by the Laplace-Adomian and
Adomian Decomposition Methods}
\author{Man Kwong Mak}
\email{mankwongmak@gmail.com}
\affiliation{Departamento de F\'{\i}sica, Facultad de Ciencias Naturales, Universidad de
Atacama, Copayapu 485, Copiap\'o, Chile}
\author{Chun Sing Leung}
\email{chun-sing-hkpu.leung@polyu.edu.hk}
\affiliation{Department of Applied Mathematics, Hong Kong Polytechnic University, Hong
Kong SAR, P. R. China}
\author{Tiberiu Harko}
\email{t.harko@ucl.ac.uk}
\affiliation{Department of Physics, Babes-Bolyai University, Kogalniceanu Street,
Cluj-Napoca 400084, Romania,}
\affiliation{School of Physics, Sun Yat-Sen University, Guangzhou 510275, People's
Republic of China}
\affiliation{Department of Mathematics, University College London, Gower Street, London
WC1E 6BT, United Kingdom}
\date{\today }

\begin{abstract}
The biharmonic equation, as well as its nonlinear and inhomogeneous
generalizations, plays an important role in engineering and physics. In
particular the focusing biharmonic nonlinear Schr\"{o}dinger equation, and
its standing wave solutions, have been intensively investigated. In the
present paper we consider the applications of the Laplace-Adomian and
Adomian Decomposition Methods for obtaining semi-analytical solutions of the
generalized biharmonic equations of the type $\Delta ^{2}y+\alpha \Delta
y+\omega y+b^{2}+g\left( y\right) =f$, where $\alpha $, $\omega $ and $b$
are constants, and $g$ and $f$ are arbitrary functions of $y$ and the
independent variable, respectively. After introducing the general algorithm
for the solution of the biharmonic equation, as an application we consider
the solutions of the one-dimensional and radially symmetric biharmonic
standing wave equation $\Delta ^{2}R+R-R^{2\sigma +1}=0$, with $\sigma =%
\mathrm{constant}$. The one-dimensional case is analyzed by using both the
Laplace-Adomian and the Adomian Decomposition Methods, respectively, and the
truncated series solutions are compared with the exact numerical solution.
The power series solution of the radial biharmonic standing wave equation is
also obtained, and compared with the numerical solution. \newline
\textit{2010 Mathematics Subject Classification}: 34K28; 34L30; 34M25;
34M30; 35C10 \newline
\textbf{Keywords:} Biharmonic equation; Laplace-Adomian Decomposition
Method; One dimensional standing wave equation; Radial standing wave equation
\end{abstract}

\pacs{02.30.Hq; 02.30.Mv; 02.30.Vv; 02.60.Cb}
\maketitle
\tableofcontents

\section{Introduction}

The biharmonic equation appears in numerous applications in science and
engineering \cite{b1b,b2b,L}. For example, the equation describing the
displacement vector $\vec{u}$ in elastodynamics is given by \cite{b2b,L}
\begin{equation}  \label{1}
\left(\lambda +\mu\right)\nabla \left(\nabla \cdot \vec{u}\right)+\mu \nabla
^2 \vec{u}+\vec{F}=0,
\end{equation}
where $\lambda $ and $\mu$ are the Lam\'{e} coefficients, and $\vec{F}$ is
the body force acting on the object. By decomposing the displacement vector $%
\vec{u}=\nabla \phi +\nabla \times \vec{\psi}$, Eq. (\ref{1}) gives
\begin{equation}
\nabla ^2\nabla ^2\phi =\nabla ^4 \phi=\Delta ^2\phi =-\frac{1}{\lambda +\mu}%
\nabla \cdot \vec{F}, \nabla ^2\nabla ^2\vec{\psi}=\nabla ^4\vec{\psi }=%
\frac{1}{\mu}\nabla \times \vec{F},
\end{equation}
that is, the equations for $\phi $ and $\vec{\psi}$ are the inhomogeneous
scalar and vector biharmonic equations \cite{b2b}.  Continuous models of elastic bodies have been intensively studied by using a variety of mathematical methods. The uniqueness of the solution of an initial-boundary value problem in thermoelasticity of bodies with voids was established in \cite{Marin1}.The theory of semigroups of operators was applied in \cite{Marin2} in order to prove the
existence and uniqueness of solutions for the mixed initial-boundary value problems
in the thermoelasticity of dipolar bodies. The temporal
behaviour of the solutions of the equations describing  a porous thermoelastic body, including voidage time
derivative among the independent constitutive variables was considered in \cite{Marin3}.

The biharmonic equation also appears in the context of gravitational
theories. Let's consider the gravitational field of Dirac $\delta$-type mass
distribution, with the mass density given by $\rho =4\pi Gm \delta \left(%
\vec{r}\right)$, where $G$ is gravitational constant, $m$ the mass, and $%
\delta \left(\vec{r}\right)$ is the Dirac delta function. Then the
gravitational potential $\Phi$ satisfies the Poisson equation \cite{Boos},
\begin{equation}
\Delta \Phi =4\pi Gm\delta \left(\vec{r}\right),
\end{equation}
with the radial solution given by $\Phi (r)=-Gm/r$. As it is well known,
this potential is singular at $r=0$, giving rise to infinite tidal forces.
However, a modification of the Poisson equation of the form \cite{Boos}
\begin{equation}  \label{gpe}
\Delta \left(1+M^{-2}\Delta\right)\Phi=4\pi Gm\delta \left(\vec{r}\right),
\end{equation}
where $M$ is a constant, gives the solution $\Phi
(r)=-Gm\left(1-e^{-Mr}\right)/r$, which is nonsingular at $r=0$, and tends
towards the Newtonian potential when $M\rightarrow \infty$.

In quantum mechanics the biharmonic equation plays an important role. The
Gross-Pitaevskii equation, describing the physical properties of
Bose-Einstein Condensates in the presence of a gravitational potential is
given by \cite{Bose1,Bose2,Bose3,Bose4}
\begin{equation}  \label{schroedinger}
i\frac{\partial}{\partial t}\,\psi \left(\vec{r}, t\right) = \left[-\frac{%
\nabla ^2}{2M^2} + \phi_{grav}\left(\vec{r}\right) + \phi_{\mathrm{rot}%
}\left(\vec{r}\right) + \phi_\eta \left(\vec{r}\right) + \frac{\partial
F(\rho)}{\partial\rho}\right]\,\psi\left(\vec{r}, t\right)),
\end{equation}
where $M$ is the mass of the particle, $\phi_{grav}$ the gravitational
potential satisfying the Poisson equation, while the potential giving the
Coriolis and centrifugal forces is given by
\begin{equation}
\phi_{\mathrm{rot}}\left(\vec{r}\right) = -\frac{1}{2}|\vec{\Omega} |^2\,|%
\vec{r}|^2 + 2\vec{\Omega} \cdot \vec{v}\times\vec{r}.
\end{equation}
The potential describing the possible viscous effects is $\phi_\eta = -\eta\;%
\vec{r}\cdot\nabla\vec{v}$ \cite{Bose5}, while $F(\rho)$ is an arbitrary
function of the particle number density, $\rho = \left|\psi\left(\vec{r}%
,t\right)\right|^2$ \cite{Bose1}. Assuming that the wave function can be
described as $\psi\left(\vec{r},t\right) = \sqrt{\rho}\,e^{i\,S\left(\vec{r}%
,t\right)}$, where $S\left(\vec{r},t\right)$ is the action of the particle,
by defining $\vec{v} = \nabla S/M $ it follows that in the static case the
Schr\"{o}dinger equation is equivalent with a system of two equations, the
continuity equation $\nabla \cdot \left(\rho \vec{v}\right)=0$, and an Euler
type equation, given by
\begin{equation}  \label{7}
\frac{1}{\rho}\,\nabla\,p + \nabla \,\left(\frac{v^2}{2} + \phi\right) +
\vec{\Omega} \times \vec{\Omega} \times \vec{r} + 2\,\vec{\Omega}\times \vec{%
v} = \eta\,\nabla ^2\,\vec{v} + \frac{1}{2M^2}\,\nabla \left(\frac{\nabla ^2%
\sqrt{\rho}}{\sqrt{\rho}}\right).
\end{equation}
This representation of the Schr\"{o}dinger equation is called the
hydrodynamic or the Madelung representation of quantum mechanics. The
pressure $p$ of the quantum fluid can be obtained from the function $F(\rho)$
as \cite{Bose1}
\begin{equation}
p = \rho\,\frac{\partial F(\rho)}{\partial \rho} - F(\rho).
\end{equation}

This relation follows from the equivalence between the Schr\"{o}dinger
equation in the hydrodynamic representation, and the Euler equation (\ref{7}%
), respectively.

In the static case, by taking the divergence of Eq.~(\ref{7}) gives a
biharmonic type equation for the density distribution of the quantum fluid,
\begin{equation}
4\pi\,G\,\rho = -\nabla \left(\frac{1}{\rho}\,\nabla \,p\right) + 2\vec{%
\Omega}^2 + \frac{1}{2M}\,\nabla ^2\left(\frac{\nabla ^2\sqrt{\rho}}{\sqrt{%
\rho}}\right).
\end{equation}

Another quantum mechanical context in which the biharmonic equation does
appear is in physical models described by the focusing biharmonic nonlinear
Schr\"{o}dinger equation, \cite{b3s,b4s,b5s,b6s,b7s,b8s},
\begin{equation}  \label{bhe}
i\hbar \frac{\partial \Psi \left(t,\vec{r}\right)}{\partial t}-\Delta ^2\Psi
\left(t,\vec{r}\right)+\left|\Psi \left(t,\vec{r}\right)\right|^{2\sigma}%
\Psi \left(t,\vec{r}\right)=0,
\end{equation}
where $\sigma \in \mathbb{R}$, and which must be solved with the initial
condition $\Psi \left(0,\vec{r}\right)=\Psi _0\left(\vec{r}\right)\in
H^2\left(\mathbb{R}^d\right)$. The focusing biharmonic nonlinear Schr\"{o}%
dinger equation is the generalization of the focusing nonlinear Schr\"{o}%
dinger equation, given by
\begin{equation}
i\hbar \frac{\partial \Psi \left(t,\vec{r}\right)}{\partial t}-\Delta \Psi
\left(t,\vec{r}\right)+\left|\Psi \left(t,\vec{r}\right)\right|^{2\sigma}%
\Psi \left(t,\vec{r}\right)=0,
\end{equation}
and it can be derived from the variational principle \cite{b3s}
\begin{equation}
S = \int \mathcal{L} d^4\vec{r} dt,
\end{equation}
where the Lagrangian density ~$\mathcal{L}$ is given by
\begin{equation}  \label{eq:Lagrangian}
\mathcal{L}\left(\psi, \psi^*, \psi_t, \psi_t^*, \Delta\psi, \Delta\psi^*
\right) = \frac{i}{2} \left( \psi_t\psi^* - \psi_t^*\psi \right) -
\left|\Delta \psi\right|^2 + \frac{1}{1+\sigma} \left|\psi\right|^{2(%
\sigma+1)}.
\end{equation}

An equation of the form
\begin{equation}
\Delta _pu+V(x)\left|u\right|^{p-2}u=f(x,u),
\end{equation}
where $p\geq 2$, and $\Delta _p^2u=\Delta \left(\left|\Delta
u\right|^{p-2}\Delta u\right)$ is called the $p$-biharmonic operator, plays
an important role in the mathematical modeling of non Newtonian fluids and
in elasticity. In particular, it describes the properties of the
electro-rheological fluids, with viscosity depending on the applied electric
field \cite{Ruz}.

Eq.~(\ref{bhe}) has the important property of admitting waveguide
(standing-wave) solutions, which can be represented as~$\psi (t,\vec{r}%
)=\lambda ^{2/\sigma }e^{i\lambda ^{4}t}R(\lambda \vec{r})$, where the
function ~$R$ satisfies the "standing-wave" equation, which takes the form
of a biharmonic equation, given by \cite{b3s}
\begin{equation}  \label{bhe1}
-\Delta ^{2}R\left( \vec{r}\right) -R\left( \vec{r}\right) +|R|^{2\sigma
}R\left( \vec{r}\right) =0.
\end{equation}%
If $\sigma d=4$, Eq.~(\ref{bhe}) is called $L^{2}$-critical, or simply
critical \cite{b3s}. The properties of the generalized nonlinear biharmonic
equation (\ref{bhe}) where studied by using mostly numerical methods \cite%
{n1,n2}. Peak-type singular solutions of Eq.~(\ref{bhe}) of the quasi-self
similar form $\Psi (t,r)\sim \left( 1/L^{d/2}(t)R\left( r/L(t)\right)
\right) e^{i\int {dt^{\prime 4}\left( t^{\prime }\right) }}$, with $%
\lim_{t\rightarrow T_{c}}L(t)=0$ have been shown to exist in \cite{b3s}.

In one dimension, Eq.~(\ref{bhe1}) is given by
\begin{equation}  \label{bhe2}
-\frac{d^4R(x)}{dx^4}-R(x)+|R|^{2\sigma}(x)R(x)=0.
\end{equation}
On the other hand, if we require radial symmetry, Eq.~(\ref{bhe1}) reduces
to
\begin{equation}  \label{bhe3}
-\Delta^2_rR(r)-R(r)+|R|^{2\sigma}(r)R(r) = 0,
\end{equation}
where~$\Delta_r^2$, the radial biharmonic operator, is given by
\begin{equation}  \label{d2}
\Delta_r^2 = \frac{d^4}{dr^4} +\frac{2(d-1)}{r}\frac{d^3}{dr^3} +\frac{%
(d-1)(d-3)}{r^2}\frac{d^2}{dr^2} -\frac{(d-1)(d-3)}{r^3}\frac{d}{dr}.
\end{equation}
At the origin $r=0$, all the odd derivatives of $R$ must vanish, and hence
the standing wave solution of the focusing biharmonic nonlinear Schr\"{o}%
dinger equation must satisfy the boundary conditions
\begin{equation}  \label{incond}
R^\prime(0)=R^{\prime\prime\prime}(0)=R(\infty)=R^\prime(\infty)=0.
\end{equation}

A lot of attention has been devoted recently to the study of Adomian's
decomposition method (ADM) \cite{new1, new2,new3,new4, R1,R2}, a powerful
mathematical method that offers the possibility of obtaining approximate
analytical solutions of many kinds of ordinary and partial differential
equations, as well as of integral equations that describe various
mathematical, physical and engineering problems. One of the important
advantages of the Adomian Decomposition Method is that it can provide
analytical approximations to the solutions of a rather large class of
nonlinear (and stochastic) differential and integral equations without the
need of linearization, or the use of perturbative and closure
approximations, or of discretization methods, which could lead to the
necessity of the extensive use of numerical computations. Usually to obtain
a closed-form analytical solutions of a nonlinear problem requires some
simplifying and restrictive assumptions.

In the case of differential equations the Adomian Decomposition Method
generates a solution in the form of a series, whose terms are obtained
recursively by using the Adomian polynomials. Together with its formal
simplicity, the main advantage of the Adomian Decomposition Method is that
the series solution of the differential equation converges fast, and
therefore its application saves a lot of computing time. Moreover, in the
Adomian Decomposition Method there is no need to discretize or linearize the
considered differential equation. For reviews of the mathematical aspects of
the Adomian Decomposition Method and its applications in physics and
engineering see \cite{R1} and \cite{R2}, respectively. From a historical point of view,
the ADM was first introduced and applied in the 1980's \cite%
{new1,new2,new3,new4}. Ever since it has been continuously modified,
generalized and extended in an attempt to improve its precision and
accuracy, and/or to expand the mathematical, physical and engineering
applications of the original method \cite%
{b2,b3,b5,b6,b7a,b8,b9,b10,b11,b12,b13,b14,b28,b30,b31,p1,p2,p3,C1,C2,
C3,C4,C5,C6,C7,C8,C9,C10}. The Adomian method was extensively applied in
mathematical physics and for the study of population growth models that can
be described by ordinary or partial differential equations, or systems of
ordinary and partial differential equations. A few example of such systems
successfully investigated by using the ADM are shallow water waves \cite%
{solw}, the Brussselator model \cite{Bruss}, the Lotka- Volterra
prey-predator type model \cite{Lotka}, and the Belousov - Zhabotinski
reduction model \cite{BJ}, respectively.  The equations of motion of the massive and massless particles in the Schwarzschild geometry of general relativity by using the Laplace-Adomian Decomposition were investigated in \cite{Mak}, where  series solutions of the geodesics equation in the Schwarzschild geometry were obtained.

Despite the considerable importance of the biharmonic equation in many
applications, very little work has been devoted to its study via the Adomian
Decomposition Method. A numerical method based on the Adomian Decomposition
Method was introduced in \cite{Khal} for the approximate solution of the one
dimensional equations of the form
\begin{equation*}
\frac{d^4u(x)}{dx^4}+\alpha (x)\frac{d^2u(x)}{dx^2}+\beta (x)\frac{du}{dx}%
=f\left(u(x)\right),
\end{equation*}
where $f\left(u(x)\right)$ is an arbitrary nonlinear function. The obtained
formalism was applied to the case of the equation
\begin{equation*}
\frac{d^4u(x)}{dx^4}+\mu u(x)=0,
\end{equation*}
where $\mu $ is a constant, and it was shown that the Adomian approximation
gives a good description of the numerical solution.

It is the purpose of the present paper to consider a systematic
investigation of the applications of the Adomian Decomposition method to the
case of the nonlinear biharmonic equation. We will consider two
implementations of the Adomian Decomposition Method, the Laplace-Adomian
Decomposition Method, and the standard Adomian Decomposition Method,
respectively. We consider both the one-dimensional nonlinear biharmonic
equation of the form
\begin{equation}
\frac{d^4 y(x)}{dx^4}+\alpha \frac{d^2y}{dx^2}+\omega y(x)+b^{2}+g(y)=f(x),
\end{equation}
as well as the nonlinear biharmonic equation with radial symmetry, given by
\begin{equation}
\frac{d^{4}y(r)}{dr^{4}}+\frac{4}{r}\frac{d^{3}y(r)}{dr^{3}}+\alpha \frac{%
d^{2}y(r)}{dr^{2}}+\frac{2}{r}\alpha \frac{dy(r)}{dr}+\omega
y(r)+b^{2}+g(y(r))=f(r).
\end{equation}%
These equations are the generalization of Eq. (\ref{bhe1}), in the
one-dimensional and radially symmetric case. For the sake of generality we
have also introduced the second order derivative whose presence allows an
easy comparison between the properties of the biharmonic and harmonic
equations. We have also included a source term in the biharmonic equations.
In both cases we develop the corresponding Laplace-Adomian and Adomian
Decomposition Method algorithms. As an application of the developed methods
we obtain the Adomian type power series solutions of the biharmonic
nonlinear standing wave equations (\ref{bhe2}) and (\ref{bhe3}),
respectively. In all cases the approximate solutions are compared with the
exact numerical ones.

The present paper is organized as follows. In Section~\ref{sect2} we discuss
the application of the Laplace-Adomian Decomposition Method to the case of
the generalized nonlinear one dimensional biharmonic equation of the type $%
\frac{d^{4}y(x)}{dx^{4}}+\alpha \frac{d^{2}y}{dx^{2}}+\omega
y(x)+b^{2}+g(y)=f(x)$. The general Laplace-Adomian Decomposition Method
algorithm is developed for this equations. As an application of our general
results we consider the one dimensional biharmonic standing wave equation $%
\frac{d^4R}{dx^{4}}+R-R^{2}=0$, and we obtain its truncated power series
solution by using both the Laplace-Adomian and the Adomian Decomposition
Methods. The truncated series solutions are compared with the exact
numerical solution. The generalized nonlinear biharmonic equation with
radial symmetry is considered in Section~\ref{sect3}. The Laplace-Adomian
Decomposition Method algorithm is developed for this case, and the solutions
of the biharmonic standing wave equation are obtained in the form of a
truncated power series. The comparison with the exact numerical solution is
also performed. Finally, we discuss and conclude our results in Section~\ref%
{sect4}.

\section{The Laplace-Adomian and the Adomian Decomposition Methods for the
nonlinear one dimensional biharmonic equation}

\label{sect2}

In the present Section we develop the Laplace-Adomian Decomposition Method
for a generalized one dimensional nonlinear inhomogeneous biharmonic type
equation of the form
\begin{equation}  \label{eqb}
\frac{d^4 y(x)}{dx^4}+\alpha \frac{d^2y}{dx^2}+\omega y(x)+b^{2}+g(y)=f(x),
\end{equation}%
where $\alpha $, $\omega $ and $b$ are constants, $g$ is an arbitrary
nonlinear function of dependent variable $y$, while $f(x)$ is an arbitrary
function of the independent variable $x$. Eq.~(\ref{eqb}) must be integrated
with the initial conditions $y(0)=y_{0}$, $y^{\prime }(0)=y_{01}$, $%
y^{\prime \prime }(0)=y_{02}$, and $y^{\prime \prime \prime}(0)=y_{03}$,
respectively.

\subsection{The general algorithm}

In the Laplace-Adomian method we apply the Laplace transformation operator $%
\mathcal{L}$, defined as $\mathcal{L}[f(x)]=\int_{0}^{\infty }f(x)e^{-sx}dx$ \cite{Lern}
, to Eq.~(\ref{eqb}). Thus we obtain
\begin{equation}
\mathcal{L}\left[ \frac{d^{4}y(x)}{dx^{4}}\right] +\alpha \mathcal{L}\left[
\frac{d^{2}y}{dx^{2}}\right] +\omega \mathcal{L}[y]+\mathcal{L}[b^{2}]+%
\mathcal{L}\left[ g(y)\right] =\mathcal{L}\left[ f(x)\right] .
\end{equation}

In the following we denote $\mathcal{L}[f(x)]=F(s)$. We use now the
properties of the Laplace transform, and thus we find
\begin{eqnarray}
F(s) &=&\frac{s\left\{ \left( \alpha +s^{2}\right) \left[ sy(0)+y^{\prime
}(0)\right] +sy^{\prime \prime }(0)+y^{\prime \prime \prime }(0)\right\}
-b^{2}}{s\left( s^{4}+\alpha s^{2}+\omega \right) }+\frac{1}{s^{4}+\alpha
s^{2}+\omega }\mathcal{L}[f(x)](s)-  \notag \\
&&\frac{1}{s^{4}+\alpha s^{2}+\omega }\mathcal{L}[g(y(x))](s).
\end{eqnarray}

As a next step we assume that the solution of the one dimensional biharmonic
Eq. (\ref{eqb}) can be represented in the form of an infinite series, given
by
\begin{equation}
y(x)=\sum_{n=0}^{\infty }y_{n}(x),  \label{7a}
\end{equation}%
where all the terms $y_{n}(x)$ can be computed recursively. As for the
nonlinear operator $g(y)$, it is decomposed according to
\begin{equation}
g(y)=\sum_{n=0}^{\infty }A_{n},  \label{8a}
\end{equation}%
where the $A_{n}$'s are the Adomian polynomials. They can be computed
generally from the definition \cite{R2}
\begin{equation}
A_{n}=\left. \frac{1}{n!}\frac{d^{n}}{d\epsilon ^{n}}f\left(
\sum_{i=0}^{\infty }{\epsilon ^{i}y_{i}}\right) \right\vert _{\epsilon =0}.
\end{equation}

The first five Adomian polynomials are given by the expressions,
\begin{equation}
A_{0}=f\left( y_{0}\right) ,  \label{Ad0}
\end{equation}%
\begin{equation}
A_{1}=y_{1}f^{\prime }\left( y_{0}\right) ,  \label{Ad1}
\end{equation}%
\begin{equation}
A_{2}=y_{2}f^{\prime }\left( y_{0}\right) +\frac{1}{2}y_{1}^{2}f^{\prime
\prime }\left( y_{0}\right) ,  \label{Ad2}
\end{equation}%
\begin{equation}
A_{3}=y_{3}f^{\prime }\left( y_{0}\right) +y_{1}y_{2}f^{\prime \prime
}\left( y_{0}\right) +\frac{1}{6}y_{1}^{3}f^{\prime \prime \prime }\left(
y_{0}\right) ,  \label{Ad3}
\end{equation}%
\begin{equation}
A_{4}=y_{4}f^{\prime }\left( y_{0}\right) +\left[ \frac{1}{2!}%
y_{2}^{2}+y_{1}y_{3}\right] f^{\prime \prime }\left( y_{0}\right) +\frac{1}{%
2!}y_{1}^{2}y_{2}f^{\prime \prime \prime }\left( y_{0}\right) +\frac{1}{4!}%
y_{1}^{4}f^{(\mathrm{iv})}\left( y_{0}\right) .  \label{Ad4}
\end{equation}

Substituting Eqs. (\ref{7a}) and (\ref{8a}) into Eq. (\ref{eqb}) we obtain
\begin{eqnarray}
\mathcal{L}\left[ \sum_{n=0}^{\infty }y_{n}(x)\right] &=&\frac{s\left\{
\left( \alpha +s^{2}\right) \left[ sy(0)+y^{\prime }(0)\right] +sy^{\prime
\prime }(0)+y^{\prime \prime \prime }(0)\right\} -b^{2}}{s\left(
s^{4}+\alpha s^{2}+\omega \right) }+  \notag  \label{11a} \\
&&\frac{1}{s^{4}+\alpha s^{2}+\omega }\mathcal{L}[f(x)](s)-\frac{1}{%
s^{4}+\alpha s^{2}+\omega }\mathcal{L}[\sum_{n=0}^{\infty }A_{n}].
\end{eqnarray}

Matching both sides of Eq. (\ref{11a}) yields the following iterative
algorithm for the power series solution of Eq. (\ref{eqb}),
\begin{equation}
\mathcal{L}\left[ y_{0}\right] =\frac{s\left\{ \left( \alpha +s^{2}\right) %
\left[ sy(0)+y^{\prime }(0)\right] +sy^{\prime \prime }(0)+y^{\prime \prime
\prime }(0)\right\} -b^{2}}{s\left( s^{4}+\alpha s^{2}+\omega \right) }+%
\frac{1}{s^{4}+\alpha s^{2}+\omega }\mathcal{L}[f(x)](s),  \label{12a}
\end{equation}%
\begin{equation}
\mathcal{L}\left[ y_{1}\right] =-\frac{1}{s^{4}+\alpha s^{2}+\omega }%
\mathcal{L}\left[ A_{0}\right] ,  \label{12b}
\end{equation}%
\begin{equation}
\mathcal{L}\left[ y_{2}\right] =-\frac{1}{s^{4}+\alpha s^{2}+\omega }%
\mathcal{L}\left[ A_{1}\right] ,  \label{12c}
\end{equation}%
\begin{equation*}
...
\end{equation*}%
\begin{equation}
\mathcal{L}\left[ y_{k+1}\right] =-\frac{1}{s^{4}+\alpha s^{2}+\omega }%
\mathcal{L}\left[ A_{k}\right] .  \label{12n}
\end{equation}

By applying the inverse Laplace transformation to Eq. (\ref{12a}), we obtain
the value of $y_{0}$. After substituting $y_{0}$ into Eq. (\ref{Ad0}), we
find easily the first Adomian polynomial $A_{0}$. Then we substitute $A_{0}$
into Eq. (\ref{12b}), and we compute the Laplace transform of the quantities
on the right-hand side of the equation. By applying the inverse Laplace
transformation we find the value of $y_{1}$. In a similar step by step
approach the other terms $y_{2}$, $y_{3}$, . . ., $y_{k+1}$, can be computed
recursively.

\subsection{Application: the one dimensional biharmonic standing wave
equation}

As an application of the previously developed Laplace-Adomian formalism we
consider the solutions of the standing wave equation (\ref{bhe1}), By
assuming the the function $R$ is real, and that $R\in \mathbb{R}_{+}$, the
standing waves equation takes the form
\begin{equation}
\frac{d^{4}R}{dx^{4}}=R^{2\sigma +1}-R.  \label{pp1}
\end{equation}%
We solve Eq. (\ref{pp1}) with the initial conditions $R^{\prime }\left(
0\right) =R^{\prime \prime \prime }\left( 0\right) =0$, and $R(0)\neq 0$ and
$R^{\prime \prime }\left( 0\right) \neq 0$, respectively. To solve Eq. (\ref%
{pp1}) we take its Laplace transform, thus obtaining
\begin{equation}
\mathcal{L}\left[ \frac{d^{4}R}{dx^{4}}\right] =\mathcal{L}\left[ R^{2\sigma
+1}-R\right] ,
\end{equation}

\begin{equation}
\left( s^{4}+1\right) \mathcal{L}\left[ R\right] =s^{3}R\left( 0\right)
+s^{2}R^{\prime }\left( 0\right) +sR^{\prime \prime }\left( 0\right)
+R^{\prime \prime \prime }\left( 0\right) +\mathcal{L}\left[ R^{2\sigma +1}%
\right] ,  \label{qq}
\end{equation}%
and%
\begin{equation}
\mathcal{L}\left[ R\right] =\frac{s^{3}R\left( 0\right) +sR^{\prime \prime
}\left( 0\right) }{s^{4}+1}+\frac{1}{s^{4}+1}\mathcal{L}\left[ R^{2\sigma +1}%
\right] ,  \label{aq}
\end{equation}%
respectively. Hence we immediately obtain

\begin{equation}
R\left( x\right) =\mathcal{L}^{-1}\left[ \frac{s^{3}R\left( 0\right)
+sR^{\prime \prime }\left( 0\right) }{s^{4}+1}\right] +\mathcal{L}%
^{-1}\left\{ \frac{1}{s^{4}+1}\mathcal{L}\left[ R^{2\sigma +1}\right]
\right\} .  \label{as}
\end{equation}%
Substituting
\begin{equation}
R\left( x\right) =\sum_{n=0}^{\infty }R_{n}\left( x\right) ,R^{2\sigma
+1}=\sum_{n=0}^{\infty }A_{n}\left( x\right) ,
\end{equation}%
where $A_{n}$ are the Adomian polynomials for all $n$, into Eq. (\ref{as})
yields
\begin{equation}
\sum_{n=0}^{\infty }R_{n}\left( x\right) =R_{0}\left( x\right)
+\sum_{n=0}^{\infty }R_{n+1}\left( x\right) =\mathcal{L}^{-1}\left[ \frac{%
s^{3}R\left( 0\right) +sR^{\prime \prime }\left( 0\right) }{s^{4}+1}\right] +%
\mathcal{L}^{-1}\left\{ \sum_{n=0}^{\infty }\left[ \frac{1}{s^{4}+1}\mathcal{%
L}\left( A_{n}\right) \right] \right\} .  \label{vv}
\end{equation}

For the function $R^{2\sigma +1}$ a few Adomian polynomials are \cite{b7an}
\begin{equation}
A_{0}=R_{0}^{2\sigma +1},  \label{h1}
\end{equation}%
\begin{equation}
A_{1}=\left( 2\sigma +1\right) R_{1}R_{0}^{2\sigma },  \label{h2}
\end{equation}%
\begin{equation}
A_{2}=\left( 2\sigma +1\right) R_{2}R_{0}^{2\sigma }+2\sigma \left( 2\sigma
+1\right) \frac{R_{1}^{2}}{2!}R_{0}^{2\sigma -1},  \label{h3}
\end{equation}%
\begin{equation}
A_{3}=\left( 2\sigma +1\right) R_{3}R_{0}^{2\sigma }+2\sigma \left( 2\sigma
+1\right) R_{1}R_{2}R_{0}^{2\sigma -1}+2\sigma \left( 2\sigma +1\right)
\left( 2\sigma -1\right) \frac{R_{1}^{3}}{3!}R_{0}^{2\sigma -2}.
\end{equation}%
We rewrite Eq. (\ref{vv}) in a recursive form as%
\begin{equation}
R_{0}\left( x\right) =\mathcal{L}^{-1}\left[ \frac{s^{3}R\left( 0\right)
+sR^{\prime \prime }\left( 0\right) }{s^{4}+1}\right] =R(0)\cos \left( \frac{%
x}{\sqrt{2}}\right) \cosh \left( \frac{x}{\sqrt{2}}\right) +R^{\prime \prime
}(0)\sin \left( \frac{x}{\sqrt{2}}\right) \sinh \left( \frac{x}{\sqrt{2}}%
\right) ,
\end{equation}%
\begin{equation}
R_{k+1}\left( x\right) =\mathcal{L}^{-1}\left\{ \frac{\mathcal{L}\left[ A_{k}%
\right] }{s^{4}+1}\right\} .
\end{equation}

For $k=0$ we have%
\begin{equation}
R_{1}\left( x\right) =\mathcal{L}^{-1}\left\{ \frac{\mathcal{L}\left[ A_{0}%
\right] }{s^{4}+1}\right\} =\mathcal{L}^{-1}\left\{ \frac{\mathcal{L}\left[
R_{0}^{1+2\sigma }\right] }{s^{4}+1}\right\} ,
\end{equation}

For $k=1$, we obtain%
\begin{equation}
R_{2}\left( x\right) =\mathcal{L}^{-1}\left\{ \frac{\mathcal{L}\left[ A_{1}%
\right] }{s^{4}+1}\right\} =\mathcal{L}^{-1}\left\{ \frac{\mathcal{L}\left[
\left( 2\sigma +1\right) R_{1}R_{0}^{2\sigma }\right] }{s^{4}+1}\right\} ,
\end{equation}%
For $k=2$, we find%
\begin{equation}
R_{3}\left( x\right) =\mathcal{L}^{-1}\left\{ \frac{\mathcal{L}\left[ A_{2}%
\right] }{s^{4}+1}\right\} =\mathcal{L}^{-1}\left\{ \frac{\mathcal{L}\left[
\left( 2\sigma +1\right) R_{2}R_{0}^{2\sigma }+\sigma \left( 2\sigma
+1\right) R_{1}^{2}R_{0}^{2\sigma -1}\right] }{s^{4}+1}\right\} ,
\end{equation}%
For $k=3$, we find%
\begin{eqnarray}
R_{4}\left( x\right) &=&\mathcal{L}^{-1}\left\{ \frac{\mathcal{L}\left(
A_{3}\right) }{s^{4}+1}\right\} =  \notag \\
&&\mathcal{L}^{-1}\left\{ \frac{\mathcal{L}\left[ \left( 2\sigma +1\right)
R_{3}R_{0}^{2\sigma }+2\sigma \left( 2\sigma +1\right)
R_{1}R_{2}R_{0}^{2\sigma -1}+\sigma \left( 2\sigma +1\right) \left( 2\sigma
-1\right) R_{1}^{3}R_{0}^{2\sigma -2}/3\right] }{s^{4}+1}\right\} .  \notag
\end{eqnarray}%
Hence the truncated semi-analytical solution of Eq. (\ref{pp1}) is given by
\begin{equation}
R\left( x\right) \approx R_{0}\left( x\right) +R_{1}\left( x\right)
+R_{2}\left( x\right) +R_{3}\left( x\right) +R_{4}\left( x\right) +....
\end{equation}

\subsubsection{The case $\protect\sigma =1/2$}

In order to give a specific example in the following we consider the case $%
\sigma =1/2$. Then the standing wave equation (\ref{pp1}) becomes
\begin{equation}  \label{ppn}
\frac{d^{4}R}{dx^{4}}=R^{2}-R.
\end{equation}

Hence we obtain the successive approximations to the solution as
\begin{eqnarray}
R_{1}(x) &=&\frac{1}{60}\Bigg\{3\left[ R(0)^{2}+\left( R^{\prime \prime
}(0)\right) ^{2}\right] \cos \left( \sqrt{2}x\right) +4\left[ 2\left(
R^{\prime \prime }(0)\right) ^{2}-5R(0)^{2}\right] \cos \left( \frac{x}{%
\sqrt{2}}\right) \times  \notag \\
&&\cosh \left( \frac{x}{\sqrt{2}}\right) +\cosh \left( \sqrt{2}x\right) %
\Bigg[\left( \left( R^{\prime \prime }(0)\right) ^{2}-R(0)^{2}\right) \cos
\left( \sqrt{2}x\right) +  \notag \\
&&3\left( R(0)^{2}+\left( R^{\prime \prime }(0)\right) ^{2}\right) \Bigg]%
+8R(0)R^{\prime \prime }(0)\sin \left( \frac{x}{\sqrt{2}}\right) \sinh
\left( \frac{x}{\sqrt{2}}\right) -  \notag \\
&&2R(0)R^{\prime \prime }(0)\sin \left( \sqrt{2}x\right) \sinh \left( \sqrt{2%
}x\right) +15\left[ R(0)-R^{\prime \prime }(0)\right] \left[ R(0)+R^{\prime
\prime }(0)\right] \Bigg\},
\end{eqnarray}%
\begin{eqnarray}
R_{2}(x) &=&\frac{1}{57600}\Bigg\{640R(0)^{3}\cos \left( \sqrt{2}x\right)
\cosh \left( \sqrt{2}x\right) -9600R(0)^{3}-384R(0)\left[ 5R(0)^{2}-4(R^{%
\prime \prime }(0))^{2}\right] \times  \notag \\
&&\cos \left( \sqrt{2}x\right) +(3+3i)\Bigg[-(64-64i)R(0)5R(0)^{2}-4(R^{%
\prime \prime }(0))^{2}\cosh \left( \sqrt{2}x\right) +  \notag \\
&&(10+5i)(R(0)+iR^{\prime \prime }(0))(R(0)-i((R^{\prime \prime
}(0))^{2}\cosh \left( \sqrt{-4+3i}x\right) +  \notag \\
&&5(R(0)+i(R^{\prime \prime }(0)))(R(0)-i(R^{\prime \prime }(0)))\Bigg(%
(2+i)(R(0)+iR^{\prime \prime }(0))\cosh \left( \sqrt{4-3i}x\right) -  \notag
\\
&&(1+2i)(R(0)-iR^{\prime \prime }(0))\cosh \left( \sqrt{4+3i}x\right) \Bigg)%
\Bigg]-6(R^{\prime \prime }(0))\left( \left( R^{\prime \prime }(0)\right)
^{2}-3R(0)^{2}\right) \times  \notag \\
&&\sin \left( \frac{3x}{\sqrt{2}}\right) \sinh \left( \frac{3x}{\sqrt{2}}%
\right) +128R^{\prime \prime }(0)\left( 3R(0)^{2}-2(R^{\prime \prime
}(0))^{2}\right) \sin \left( \sqrt{2}x\right) \sinh \left( \sqrt{2}x\right) +
\notag \\
&&6R(0)\left( R(0)^{2}-3(R^{\prime \prime }(0))^{2}\right) \cos \left( \frac{%
3x}{\sqrt{2}}\right) \cosh \left( \frac{3x}{\sqrt{2}}\right) +  \notag \\
&&2\cos \left( \frac{x}{\sqrt{2}}\right) \Bigg[R(0)6367R(0)^{2}-1557(R^{%
\prime \prime }(0))^{2}\cosh \left( \frac{x}{\sqrt{2}}\right) -30\sqrt{2}%
x(R(0)-R^{\prime \prime }(0))\times  \notag \\
&&71R(0)^{2}+116R(0)(R^{\prime \prime }(0))+71(R^{\prime \prime
}(0))^{2}\sinh \left( \frac{x}{\sqrt{2}}\right) \Bigg]+  \notag \\
&&2\sin \left( \frac{x}{\sqrt{2}}\right) \Bigg[R^{\prime \prime }(0)\left(
4821R(0)^{2}-7711\left(R^{\prime \prime}(0)\right)^ 2\right) \sinh \left(
\frac{x}{\sqrt{2}}\right) +30\sqrt{2}x(R(0)+R^{\prime \prime }(0))\times
\notag \\
&&\Bigg(71R(0)^{2}-116R(0)R^{\prime \prime }(0)+71(R^{\prime \prime }(0))^{2}%
\Bigg)\cosh \left( \frac{x}{\sqrt{2}}\right) \Bigg]+  \notag \\
&&(15-45i)\left( R(0)-iR^{\prime \prime }(0)\right) \left( R(0)+iR^{\prime
\prime }(0)\right) ^{2}\cosh \left( \sqrt{-4-3i}x\right) \Bigg\}.
\end{eqnarray}

Thus we have obtained the following three terms truncated approximate
solution of the nonlinear one dimensional biharmonic equation (\ref{ppn}),
\begin{equation}  \label{lsol}
R(x)\approx R_0(x)+R_1(x)+R_2(x).
\end{equation}

\subsection{The Adomian Decomposition Method for the biharmonic standing
wave equation}

For the sake of comparison we also consider the application of the standard
Adomian Decomposition Method for solving the standing wave equation (\ref%
{pp1}) with the same initial conditions as used in the previous Section.
Four fold integrating Eq.~(\ref{pp1}) gives
\begin{equation}
\int_{0}^{x}dx_{1}\int_{0}^{x_{1}}dx_{2}\int_{0}^{x_{2}}dx_{3}%
\int_{0}^{x_{3}}\frac{d^{4}R\left( x_{4}\right) }{dx_{4}^{4}}%
dx_{4}=\int_{0}^{x}dx_{1}\int_{0}^{x_{1}}dx_{2}\int_{0}^{x_{2}}dx_{3}%
\int_{0}^{x_{3}}\left[ R^{2\sigma +1}\left( x_{4}\right) -R\left(
x_{4}\right) \right] dx_{4}.  \label{p5}
\end{equation}%
Hence we immediately obtain
\begin{equation}
R\left( x\right) =R\left( 0\right) +R^{\prime \prime }\left( 0\right) \frac{%
x^{2}}{2}+\int_{0}^{x}dx_{1}\int_{0}^{x_{1}}dx_{2}\int_{0}^{x_{2}}dx_{3}%
\int_{0}^{x_{3}}\left[ R^{2\sigma +1}\left( x_{4}\right) -R\left(
x_{4}\right) \right] dx_{4}.  \label{bb}
\end{equation}%
Substituting $R\left( x\right) =\sum_{n=0}^{\infty }R_{n}\left( x\right) $, $%
R^{2\sigma +1}=\sum_{n=0}^{\infty }A_{n}\left( x\right) $, where $A_{n}$ are
the Adomian polynomials for all $n$, into Eq. (\ref{bb}), yields
\begin{eqnarray}
\sum_{n=0}^{\infty }R_{n}\left( x\right) &=&R_{0}\left( x\right)
+\sum_{n=0}^{\infty }R_{n+1}\left( x\right) =R\left( 0\right) +R^{\prime
\prime }\left( 0\right) \frac{x^{2}}{2}+  \notag  \label{sa} \\
&&\sum_{n=0}^{\infty
}\int_{0}^{x}dx_{1}\int_{0}^{x_{1}}dx_{2}\int_{0}^{x_{2}}dx_{3}%
\int_{0}^{x_{3}}\left[ A_{n}\left( x_{4}\right) -R_{n}\left( x_{4}\right) %
\right] dx_{4}.
\end{eqnarray}%
We rewrite Eq. (\ref{sa}) in a recursive form as
\begin{equation}
R_{0}\left( x\right) =R\left( 0\right) +R^{\prime \prime }\left( 0\right)
\frac{x^{2}}{2},  \label{a1}
\end{equation}%
\begin{equation}
R_{k+1}\left( x\right)
=\int_{0}^{x}dx_{1}\int_{0}^{x_{1}}dx_{2}\int_{0}^{x_{2}}dx_{3}%
\int_{0}^{x_{3}}\left[ A_{k}\left( x_{4}\right) -R_{k}\left( x_{4}\right) %
\right] dx_{4}.  \label{a2}
\end{equation}

With the help of Eq.~(\ref{a1}) and Eq.~(\ref{a2}), we obtain the
semi-analytical solution of Eq. (\ref{pp1}) as given by
\begin{equation}
R\left( x\right) =R_{0}\left( x\right) +R_{1}\left( x\right) +R_{2}\left(
x\right) +R_{3}\left( x\right) ....
\end{equation}

In order to discuss a specific case we consider again Eq.~(\ref{pp1}) for $%
\sigma =1/2$. Then
\begin{equation}
A_{0}=R_{0}^{2}=\left[ R\left( 0\right) +R^{\prime \prime }\left( 0\right)
\frac{x^{2}}{2}\right] ^{2},
\end{equation}%
which gives
\begin{eqnarray}
R_{1}(x)
&=&\int_{0}^{x}dx_{1}\int_{0}^{x_{1}}dx_{2}\int_{0}^{x_{2}}dx_{3}%
\int_{0}^{x_{3}}\left[ A_{0}\left( x_{4}\right) -R_{0}\left( x_{4}\right) %
\right] dx_{4}= \\
&&\int_{0}^{x}dx_{1}\int_{0}^{x_{1}}dx_{2}\int_{0}^{x_{2}}dx_{3}%
\int_{0}^{x_{3}}\left\{ \left[ R\left( 0\right) +R^{\prime \prime }\left(
0\right) \frac{x_{4}^{2}}{2}\right] ^{2}-\left[ R\left( 0\right) +R^{\prime
\prime }\left( 0\right) \frac{x_{4}^{2}}{2}\right] \right\} dx_{4},
\end{eqnarray}%
\begin{equation}
R_{1}(x)=\frac{1}{24}\left[ R(0)-1\right] R(0)x^{4}+\frac{1}{720}\left[
2R(0)-1\right] R^{\prime \prime 6}+\frac{\left( R^{\prime \prime }(0)\right)
^{2}x^{8}}{6720},
\end{equation}

\begin{equation}
R_{2}\left( x\right)
=\int_{0}^{x}dx_{1}\int_{0}^{x_{1}}dx_{2}\int_{0}^{x_{2}}dx_{3}%
\int_{0}^{x_{3}}\left[ 2R_{0}\left( x_{4}\right) R_{1}\left( x_{4}\right)
-R_{1}\left( x_{4}\right) \right] dx_{4},
\end{equation}%
\begin{eqnarray}
R_{2}\left( x\right) &=&\frac{R(0)\left[ 2R(0)^{2}-3R(0)+1\right] }{40320}%
x^{8}+\frac{\left[ 34R(0)^{2}-34R(0)+1\right] R^{\prime \prime }(0)}{3628800}%
x^{10}+  \notag \\
&&\frac{31\left[ 2R(0)-1\right] \left( R^{\prime \prime }(0)\right) ^{2}}{%
239500800}x^{12}+\frac{\left( R^{\prime \prime }(0)\right) ^{3}}{161441280}%
x^{14},
\end{eqnarray}%
\begin{equation}
R_{3}\left( x\right)
=\int_{0}^{x}dx_{1}\int_{0}^{x_{1}}dx_{2}\int_{0}^{x_{2}}dx_{3}%
\int_{0}^{x_{3}}\left[ 2R_{0}\left( x_{4}\right) R_{2}\left( x_{4}\right)
+R_{1}^{2}\left( x_{4}\right) -R_{2}\left( x_{4}\right) \right] dx_{4},
\end{equation}%
\begin{eqnarray}
R_{3}(x) &=&\frac{R(0)\left[ 74R(0)^{3}-148R(0)^{2}+75R(0)-1\right] }{%
479001600}x^{12}+\frac{\left[ 1088R(0)^{3}-1632R(0)^{2}+546R(0)-1\right]
R^{\prime \prime }(0)}{87178291200}x^{14}+  \notag \\
&&\frac{\left[ 7186R(0)^{2}-7186R(0)+559\right] \left( R^{\prime \prime
}(0)\right) ^{2}}{10461394944000}x^{16}+\frac{2393\left[ 2R(0)-1\right]
\left( R^{\prime \prime }(0)\right) ^{3}}{320118685286400}x^{18}+\frac{%
61\left( R^{\prime \prime }(0)\right) ^{4}}{250298560512000}x^{20}.
\end{eqnarray}

Thus we have obtained an approximate solution of Eq.~(\ref{ppn}) as given by
\begin{equation}  \label{asol}
R(x)\approx R_0(x)+R_1(x)+R_2(x)+R_3(x).
\end{equation}

\subsection{Comparison with the exact numerical solution}

\begin{figure*}[htp!]
\centering
\includegraphics[width=8.5cm]{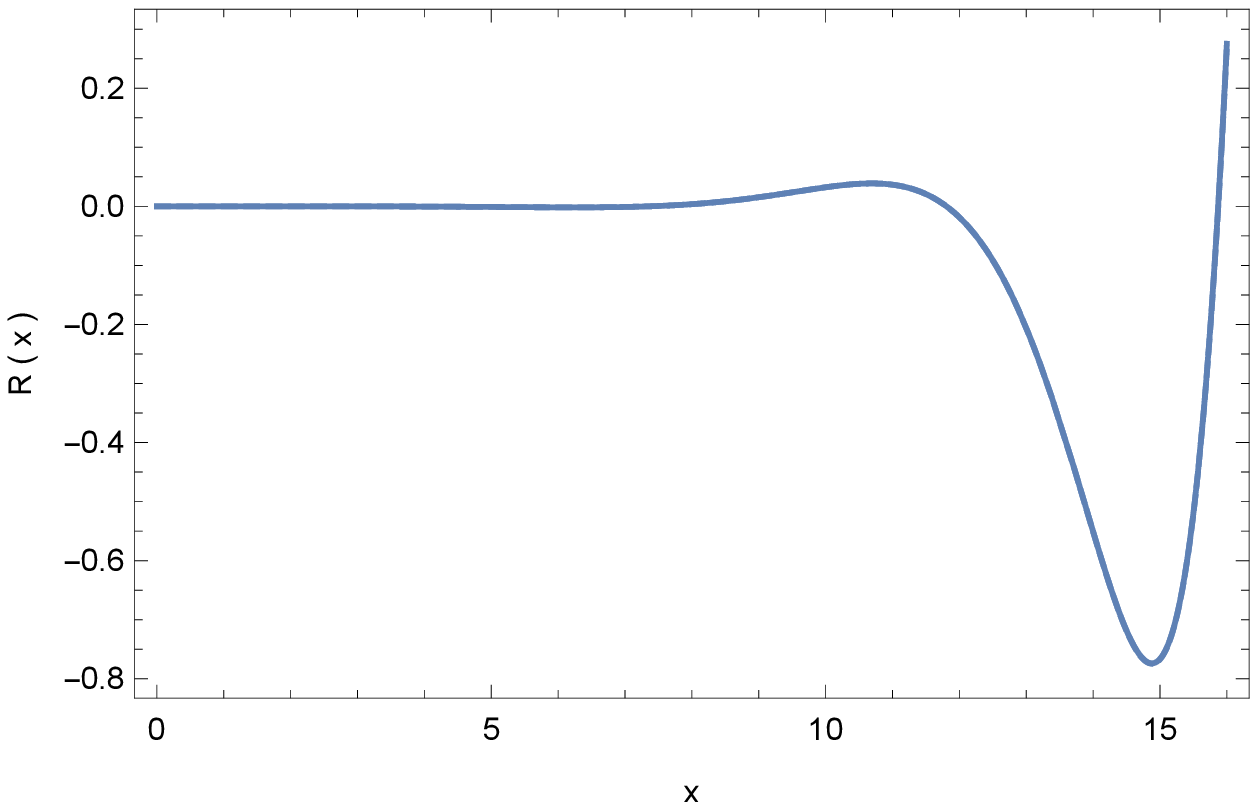} %
\includegraphics[width=8.5cm]{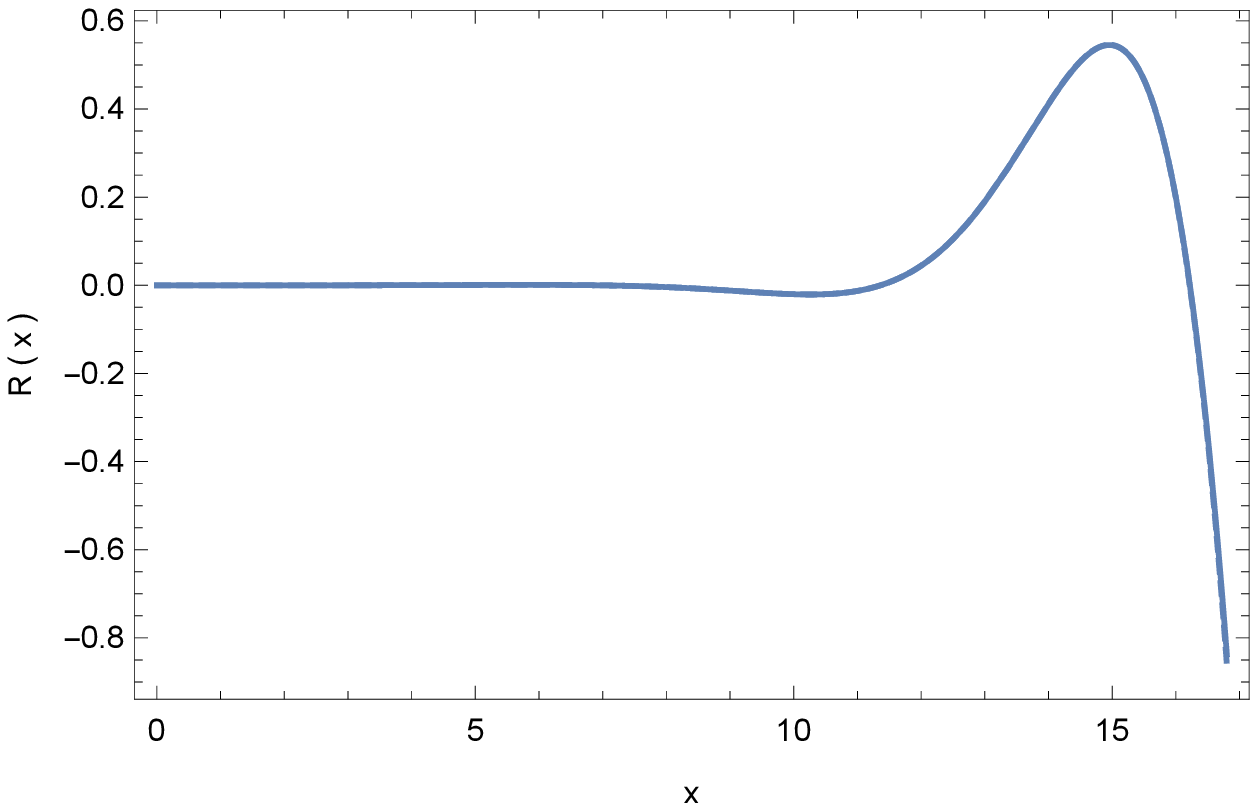}
\caption{Comparison of the numerical solutions of the nonlinear biharmonic
standing wave equation (\protect\ref{ppn}) and of the Laplace-Adomian
Decomposition Method approximate solutions, truncated to three terms, given
by Eq.~(\protect\ref{lsol}. The numerical solution is represented by the
solid curve, while the dashed curve depicts the Laplace-Adomian three terms
solution. The initial conditions used to integrate the equations are $%
R(0)=5.1\times 10^{-5}$ and $R^{\prime \prime}(0)=2.65\times 10^{ -5}$ (left
figure), and $R(0)=-4.1\times 10^{-5}$ and $R^{\prime \prime}(0)=-7.86\times
10^{-6}$ (right figure), respectively.}
\label{fig1}
\end{figure*}

\begin{figure*}[htp!]
\centering
\includegraphics[width=8.5cm]{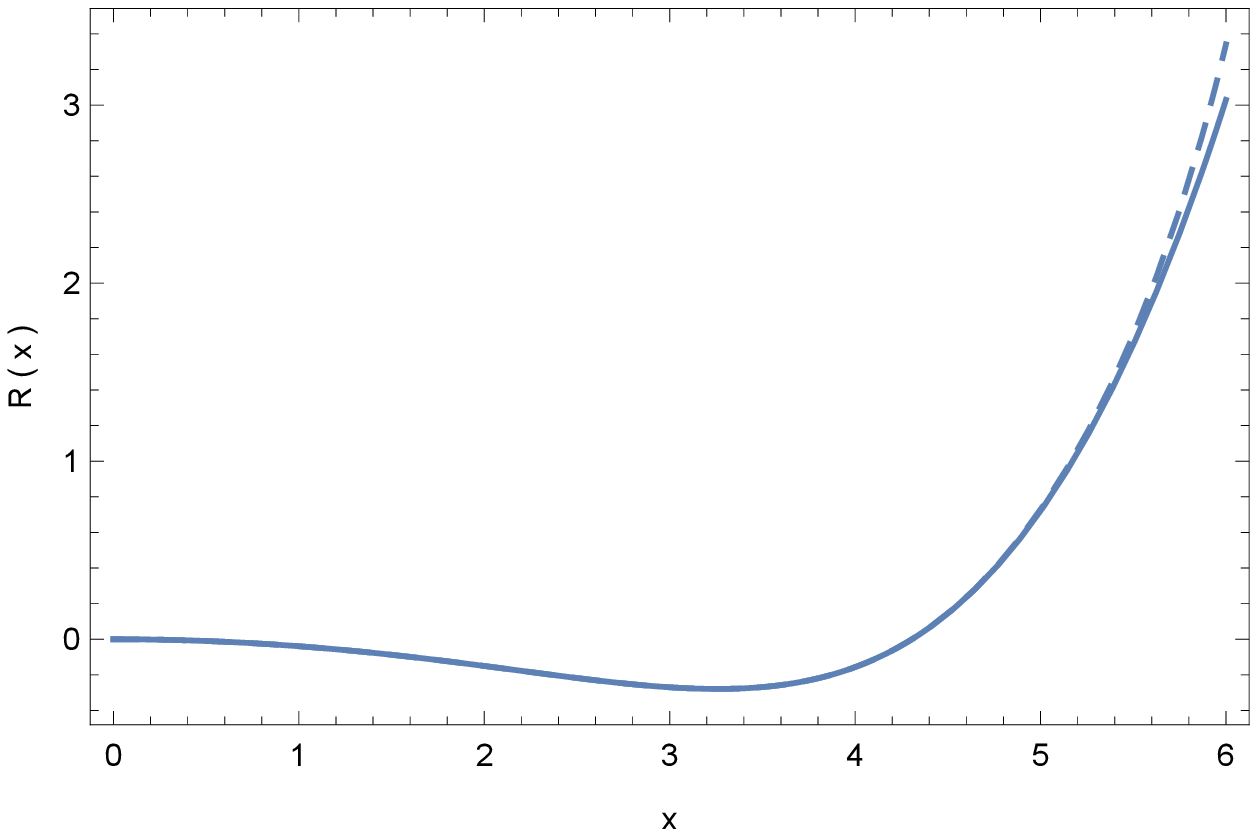} %
\includegraphics[width=8.5cm]{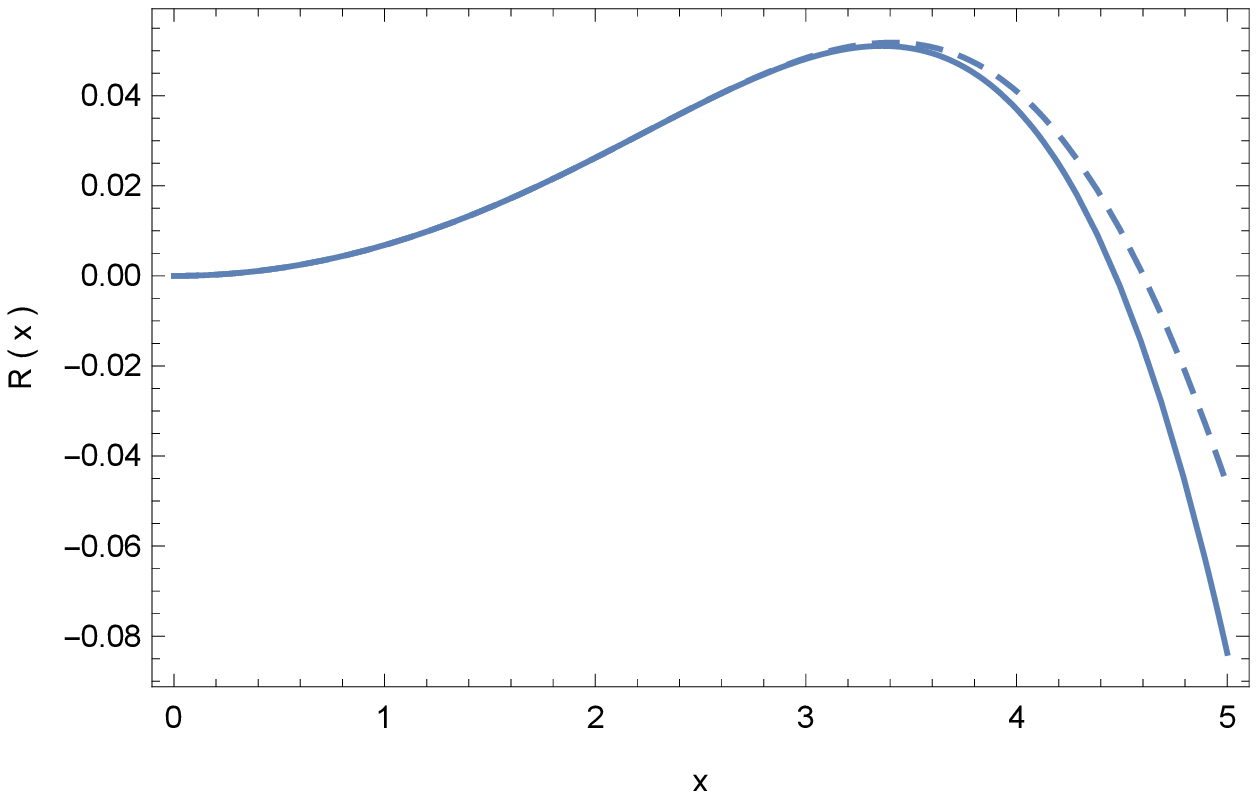}
\caption{Comparison of the numerical solutions of the nonlinear biharmonic
equation (\protect\ref{ppn}) and of the Adomian Decomposition Method
approximate solutions, truncated to four terms, given by Eq.~(\protect\ref%
{asol}. The numerical solutions are represented by the solid curves, while
the dashed curves depicts the Adomian Decomposition Method four terms
solutions. The initial conditions used to integrate the equations are $%
R(0)=-4.1\times 10^{-6}$ and $R^{\prime \prime}(0)=-7.86\times 10^{-2}$
(left figure) and $R(0)=7.19\times 10^{-8}$ and $R^{\prime
\prime}(0)=1.37\times 10^{ -2}$ (right figure), respectively.}
\label{fig2}
\end{figure*}

In order to test the accuracy of the obtained semi-analytical solutions of
the standing wave equation (\ref{pp1}) we compare the exact numerical
solution of the equation for $\sigma =1/2$ with the approximate solutions
obtained via the Laplace-Adomian and Adomian Decomposition Method. The
comparison of the exact numerical solution and the three-terms solution of
the Laplace-Adomian Method is presented in Fig.~\ref{fig1}, while the
comparison of the numerical solution and Adomian Decomposition Method is
done in Fig.~\ref{fig2}.

As one can see from Fig~\ref{fig1}, the Laplace Adomian Decomposition
Method, truncated to three terms only, gives an excellent description of the
numerical solution, at least for the adopted range of initial conditions.
The approximate solutions describes well the complex features of the
solution on a relatively large range of the independent variable $x$. The
simple Adomian Decomposition Method is more easy to apply, however, its
accuracy seems to be limited, as compared to the Laplace Adomian
Decomposition Method. Moreover, it is important to point out that there is a
strong dependence on the initial conditions of the accuracy of the method.
If the values $R(0)$ and $R^{\prime \prime }(0)$ are small, the series
solutions are in good agreement with the numerical ones. However, for larger
values of the initial conditions, the accuracy of the Adomian methods
decreases rapidly.

\section{The biharmonic nonlinear equation with radial symmetry}

\label{sect3}

In three dimensions $d=3$, and the radial biharmonic operator (\ref{d2})
takes the simple form
\begin{equation}
\Delta _r^2=\frac{d^4}{dr^4}+\frac{4}{r}\frac{d^3}{dr^3}.
\end{equation}

Hence the general nonlinear three dimensional biharmonic equation with
radial symmetry is given by
\begin{equation}
\frac{d^{4}y(r)}{dr^{4}}+\frac{4}{r}\frac{d^{3}y(r)}{dr^{3}}+\alpha \frac{%
d^{2}y(r)}{dr^{2}}+\frac{2}{r}\alpha \frac{dy(r)}{dr}+\omega
y(r)+b^{2}+g(y(r))=f(r),  \label{rad}
\end{equation}%
where $\alpha $, $b^{2}$ and $\omega $ are constants, while $g(y)$, the
nonlinear operator term, and $f(r)$, are two arbitrary functions. Eq.~(\ref%
{rad}) must be integrated with the initial conditions $y(0)=y_{0}$, $%
y^{\prime }(0)=y_{01}$, $y^{\prime \prime }(0)=y_{02}$, and $y^{\prime
\prime \prime }(0)=y_{03}$, respectively. After multiplying Eq. (\ref{rad})
with $r$ we obtain

\begin{equation}
r\frac{d^{4}y(r)}{dr^{4}}+4\frac{d^{3}y(r)}{dr^{3}}+\alpha r\frac{d^{2}y(r)}{%
dr^{2}}+2\alpha \frac{dy(r)}{dr}+\omega ry(r)+rg(y(r))=rf(r)-b^{2}r.
\label{rad1}
\end{equation}

\subsection{The Laplace-Adomian Decomposition Method solution}

As a first step in our study we assume that $y$ and $g(y(r))$ can be
represented in the form of a power series as

\begin{equation}
y=\sum_{n=0}^{\infty }y_{n},g(y(r))=\sum_{n=0}^{\infty }A_{n},
\end{equation}%
where $A_{n}$ are the Adomian polynomials. Hence Eq. (\ref{rad1}) becomes
\begin{eqnarray}  \label{rad1a}
&&\sum_{n=0}^{\infty }r\frac{d^{4}y_{n}(r)}{dr^{4}}+4\sum_{n=0}^{\infty }%
\frac{d^{3}y_{n}(r)}{dr^{3}}+\alpha \sum_{n=0}^{\infty }r\frac{d^{2}y_{n}(r)%
}{dr^{2}}+2\alpha \sum_{n=0}^{\infty }\frac{dy_{n}(r)}{dr}+  \notag \\
&&\omega \sum_{n=0}^{\infty }ry_{n}(r)+\sum_{n=0}^{\infty
}rA_{n}=rf(r)-b^{2}r.
\end{eqnarray}

After applying the Laplace transformation operator to Eq. (\ref{rad1}) we
obtain
\begin{eqnarray}
&&\sum_{n=0}^{\infty }\mathcal{L}\left[ r\frac{d^{4}y_{n}(r)}{dr^{4}}\right]
+4\sum_{n=0}^{\infty }\mathcal{L}\left[ \frac{d^{3}y_{n}(r)}{dr^{3}}\right]
+\alpha \sum_{n=0}^{\infty }\mathcal{L}\left[ r\frac{d^{2}y_{n}(r)}{dr^{2}}%
\right] +2\alpha \sum_{n=0}^{\infty }\mathcal{L}\left[ \frac{dy_{n}(r)}{dr}%
\right] +  \notag  \label{rad2} \\
&&\omega \sum_{n=0}^{\infty }\mathcal{L}\left[ ry_{n}(r)\right]
+\sum_{n=0}^{\infty }\mathcal{L}\left[ rA_{n}\right] =\mathcal{L}\left[
rf(r)-b^{2}r\right] .
\end{eqnarray}

By taking into account the relations

\begin{equation}
\mathcal{L}\left[ r\frac{d^{4}y_{n}(r)}{dr^{4}}\right] (s)=\int_{0}^{\infty
}r\frac{d^{4}y_{n}(r)}{dr^{4}}e^{-sr}dr=-\frac{d}{ds}\int_{0}^{\infty }\frac{%
d^{4}y_{n}(r)}{dr^{4}}e^{-sr}dr=-\frac{d}{ds}\mathcal{L}\left[ \frac{%
d^{4}y_{n}(r)}{dr^{4}}\right] (s),
\end{equation}

\begin{equation}
\mathcal{L}\left[ r\frac{d^{2}y_{n}(r)}{dr^{2}}\right] =\int_{0}^{\infty }r%
\frac{d^{2}y_{n}(r)}{dr^{2}}e^{-sr}dr=-\frac{d}{ds}\int_{0}^{\infty }\frac{%
d^{2}y_{n}(r)}{dr^{2}}e^{-sr}dr=-\frac{d}{ds}\mathcal{L}\left[ \frac{%
d^{2}y_{n}(r)}{dr^{2}}\right] (s),
\end{equation}

\begin{equation}
\mathcal{L}\left[ ry_{n}(r)\right] (s)=\int_{0}^{\infty }ry_{n}(r)e^{-sr}dr=-%
\frac{d}{ds}\int_{0}^{\infty }y_{n}(r)e^{-sr}dr=-\frac{d}{ds}\mathcal{L}%
\left[ y_{n}[(r)\right] (s),
\end{equation}%
and the linearity of the Laplace transformation, Eq. (\ref{rad2}) becomes
\begin{eqnarray}
&&-\left( s^{4}+\alpha s^{2}+\omega \right) \sum_{n=0}^{\infty
}F_{n}^{\prime }(s)-y(0)\left( \alpha +s^{2}\right) -2sy^{\prime
}(0)-3y^{\prime \prime }(0)+  \notag  \label{rad3} \\
&&\sum_{n=0}^{\infty }\mathcal{L}\left[ rA_{n}\right] (s)=-\frac{b^{2}}{s^{2}%
}+\mathcal{L}\left[ rf(r)\right] (s).
\end{eqnarray}

From Eq. (\ref{rad3}) we obtain the following recursion relations

\begin{equation}
-\left( s^{4}+\alpha s^{2}+\omega \right) F_{0}^{\prime }(s)-y(0)\left(
\alpha +s^{2}\right) -2sy^{\prime }(0)-3y^{\prime \prime }(0)=-\frac{b^{2}}{%
s^2}+\mathcal{L}\left[ rf(r)\right] (s),  \label{rad5}
\end{equation}

\begin{equation}
F_{n+1}^{\prime }(s)=\frac{1}{\left( s^{4}+\alpha s^{2}+\omega \right) }%
\mathcal{L}\left[ rA_{n}\right] (s).  \label{rad6}
\end{equation}

From Eq. (\ref{rad5}) we obtain
\begin{equation}
F_{0}(s)=\int G(s)ds,
\end{equation}%
where
\begin{equation}
G(s)=\frac{b^{2}/s^{2}-y(0)\left( \alpha +s^{2}\right) -2y^{\prime
}(0)s-3y^{\prime \prime }(0)-\mathcal{L}\left[ rf(r)\right] (s)}{\left(
s^{4}+\alpha s^{2}+\omega \right) },
\end{equation}%
while Eq. (\ref{rad6}) gives
\begin{equation}
F_{k+1}(s)=\int \frac{1}{\left( s^{4}+\alpha s^{2}+\omega \right) }\mathcal{L%
}\left[ rA_{k}\right] (s)ds.
\end{equation}

Hence we obtain the following approximate series solution of the radial
nonlinear biharmonic equation (\ref{rad}),
\begin{equation}
y_{0}(r)=\mathcal{L}^{-1}\left[ \int G(s)ds\right] (r),
\end{equation}
\begin{equation}
y_{k+1}=\mathcal{L}^{-1}\left\{ \int \frac{1}{\left( s^{4}+\alpha
s^{2}+\omega \right) }\mathcal{L}\left[ rA_{k}\right] (s)ds\right\} (r).
\end{equation}

\subsection{Application: the radial biharmonic standing wave equation}

In radial symmetry, and by assuming that $R\in \mathbb{R}_{+}$, the standing
wave equation (\ref{bhe3}) takes the form
\begin{equation}  \label{sw}
\frac{d^{4}R(r)}{dr^{4}}+\frac{4}{r}\frac{d^{3}R(r)}{dr^{3}}+R(r)-R^{2\sigma
+1}(r)=0,
\end{equation}%
or, equivalently,
\begin{equation}
r\frac{d^{4}R(r)}{dr^{4}}+4\frac{d^{3}R(r)}{dr^{3}}+rR(r)=rR^{2\sigma +1}(r).
\label{rad7}
\end{equation}%
By taking the Laplace transform of Eq.~(\ref{rad7}) we obtain
\begin{equation}
-\left( s^{4}+1\right) F^{\prime }(s)-R(0)s^{2}-2R^{\prime }(0)s-3R^{\prime
\prime }(0)=\mathcal{L}\left[ rR^{2\sigma +1}(r)\right] (s).  \label{rad8}
\end{equation}

By writing $R(r)=\sum_{n=0}^{\infty }R_{n}(r)$, $\mathcal{L}\left[ R(r)%
\right] (s)=\sum_{n=0}^{\infty }\mathcal{L}\left[ R_{n}(r)\right]
(s)=\sum_{n=0}^{\infty }F_{n}(s)$, $R^{2\sigma +1}(r)=\sum_{n=0}^{\infty
}A_{n}(r)$, Eq.~(\ref{rad8}) becomes
\begin{equation}
-\left( s^{4}+1\right) \sum_{n=0}^{\infty }F_{n}^{\prime
}(s)-R(0)s^{2}-2R^{\prime }(0)s-3R^{\prime \prime }(0)=\sum_{n=0}^{\infty }%
\mathcal{L}\left[ rA_{n}(r)\right] (s).
\end{equation}

Hence we obtain the following recursive relations for the solution of Eq. (%
\ref{sw}),
\begin{equation}
F_{0}^{\prime }(s)=-\frac{R(0)s^{2}+2R^{\prime }(0)s+3R^{\prime \prime }(0)}{%
s^{4}+1},  \label{F0s}
\end{equation}%
\begin{equation}
F_{k+1}^{\prime }(s)=-\frac{1}{s^{4}+1}\mathcal{L}\left[ rA_{k}(r)\right]
(s).
\end{equation}

Eq.~(\ref{F0s}) can be integrated exactly to obtain $F_{0}(s)$ as
\begin{eqnarray}  \label{91}
F_{0}(s) &=&\frac{1}{8}\Bigg\{2\left[ \sqrt{2}R(0)+4R^{\prime }(0)+3\sqrt{2}%
R^{\prime \prime }(0)\right] \tan ^{-1}\left( 1-\sqrt{2}s\right) -2\left[
\sqrt{2}R(0)-4R^{\prime }(0)+3\sqrt{2}R^{\prime \prime }(0)\right] \times
\notag \\
&&\tan ^{-1}\left( 1+\sqrt{2}s\right) - \sqrt{2}\left[ R(0)-3R^{\prime
\prime }(0)\right] \ln \frac{s^{2}-\sqrt{2}s+1}{s^{2}+\sqrt{2}s+1}\Bigg\}.
\end{eqnarray}

In the following we will consider the solutions of Eq.~(\ref{sw}) with $%
\sigma =1/2$, together with the initial conditions $R(0)\neq 0$, $R^{\prime
}(0)=0$, $R^{\prime \prime }(0)\neq 0$, and $R^{\prime \prime \prime }(0)=0$%
, respectively. Then, by neglecting the non-linear term $R^{2}$ in Eq.~(\ref%
{sw}) it turns out that the general solution of the linear equation
\begin{equation}
r\frac{d^{4}R_{0}(r)}{dr^{4}}+4\frac{d^{3}R_{0}(r)}{dr^{3}}+rR_{0}(r)=0,
\end{equation}%
is given by
\begin{equation}
R_{0}(r)=\frac{\left( \frac{1}{2}-\frac{i}{2}\right) \left\{ \sin \left(
\sqrt[4]{-1}r\right) \left[ R(0)+3iR^{\prime \prime }(0)\right] +\sinh
\left( \sqrt[4]{-1}r\right) \left[ R(0)-3iR^{\prime \prime }(0)\right]
\right\} }{\sqrt{2}r}.
\end{equation}%
The Laplace transform of $R_{0}$ converges only for values of $\mathrm{Re}%
\;s\geq s_{0}=1/\sqrt{2}$. In the region of convergence $F_{0}(s)$ can
effectively be expressed as the absolutely convergent Laplace transform of
another function, such that $F_{0}(s)=\left( s-s_{0}\right) \int_{0}^{\infty
}{e^{-\left( s-s_{0}\right) t}\beta (t)dt}$, where $\beta (t)=\int_{0}^{t}{%
e^{-s_{0}u}R_{0}(u)du}$.

The first Adomian polynomial $A_0$ is obtained as $A_0(r)=R_0^2(r)$, and the
Laplace transform of $rA_0$ is given by
\begin{eqnarray}
\mathcal{L}\left[ rA_{0}(r)\right](s)&=&\frac{1}{8} \Bigg\{3 R(0) R^{\prime
\prime }(0) \log \left(\frac{16}{s^4}+1\right)+ \left[R(0)^2+9
\left(R^{\prime \prime }(0)\right)^2\right]\ln \left(\frac{s^2+2}{s^2-2}%
\right)+  \notag \\
&& \left[R(0)^2-9 \left(R^{\prime \prime }(0)\right)^2\right]\tan ^{-1}\left(%
\frac{4}{s^2}\right) \Bigg\}.
\end{eqnarray}

Then the Laplace transform of the first correction term in the Adomian
series expansion is given as the solution of the following differential
equation,
\begin{eqnarray}
F_{1}^{\prime }(s) &=&-\frac{1}{8\left( 1+s^{4}\right) }\Bigg\{%
3R(0)R^{\prime \prime }(0)\log \left( \frac{16}{s^{4}}+1\right) +\left[
R(0)^{2}+9\left( R^{\prime \prime }(0)\right) ^{2}\right] \ln \left( \frac{%
s^{2}+2}{s^{2}-2}\right) +  \notag \\
&&\left[ R(0)^{2}-9\left( R^{\prime \prime }(0)\right) ^{2}\right] \tan
^{-1}\left( \frac{4}{s^{2}}\right) \Bigg\}.
\end{eqnarray}%
The right hand side of the above equation cannot be integrated exactly. By
expanding it in power series of $1/s$, we obtain
\begin{eqnarray}
F_{1}^{\prime }(s) &\approx &-\frac{R(0)^{2}}{s^{6}}-\frac{6(R(0)R^{\prime
\prime }(0))}{s^{8}}+\frac{3R(0)^{2}-30\left( R^{\prime \prime }(0)\right)
^{2}}{s^{10}}+\frac{54R(0)R^{\prime \prime }(0)}{s^{12}}+\frac{246\left(
R^{\prime \prime }(0)\right) ^{2}-\frac{151R(0)^{2}}{5}}{s^{14}}-  \notag \\
&&\frac{566(R(0)R^{\prime \prime }(0))}{s^{16}}+O\left( \left( \frac{1}{s}%
\right) ^{17}\right) ,
\end{eqnarray}%
and
\begin{eqnarray}
F_{1}(s) &\approx &\frac{1}{5}\Bigg\{\frac{566R(0)R^{\prime \prime }(0)}{%
3s^{15}}+\frac{151R(0)^{2}-1230\left( R^{\prime \prime }(0)\right) ^{2}}{%
13s^{13}}-\frac{270R(0)R^{\prime \prime }(0)}{11s^{11}}-\frac{5\left[
R(0)^{2}-10\left( R^{\prime \prime }(0)\right) ^{2}\right] }{3s^{9}}+  \notag
\\
&&\frac{30R(0)R^{\prime \prime }(0)}{7s^{7}}+\frac{R(0)^{2}}{s^{5}}\Bigg\}.
\end{eqnarray}%
respectively. Hence for the first term of the Adomian series expansion $%
R_{1}=\mathcal{L}^{-1}\left[ F_{1}(s)\right] (r)$ we obtain
\begin{eqnarray}
R_{1}(r) &=&\frac{R(0)^{2}}{120}r^{4}+\frac{1}{840}R(0)R^{\prime \prime
}(0)r^{6}+\frac{ 10\left( R^{\prime \prime }(0)\right) ^{2}-R(0)^{2} }{120960%
}r^{8}-\frac{R(0)R^{\prime \prime }(0)}{739200}r^{10}+  \notag \\
&&\frac{ 151R(0)^{2}-1230\left( R^{\prime \prime }(0)\right) ^{2} }{%
31135104000}r^{12}+\frac{283R(0)R^{\prime \prime }(0)}{653837184000}%
r^{14}+....
\end{eqnarray}

\begin{figure*}[htp!]
\centering
\includegraphics[width=8.5cm]{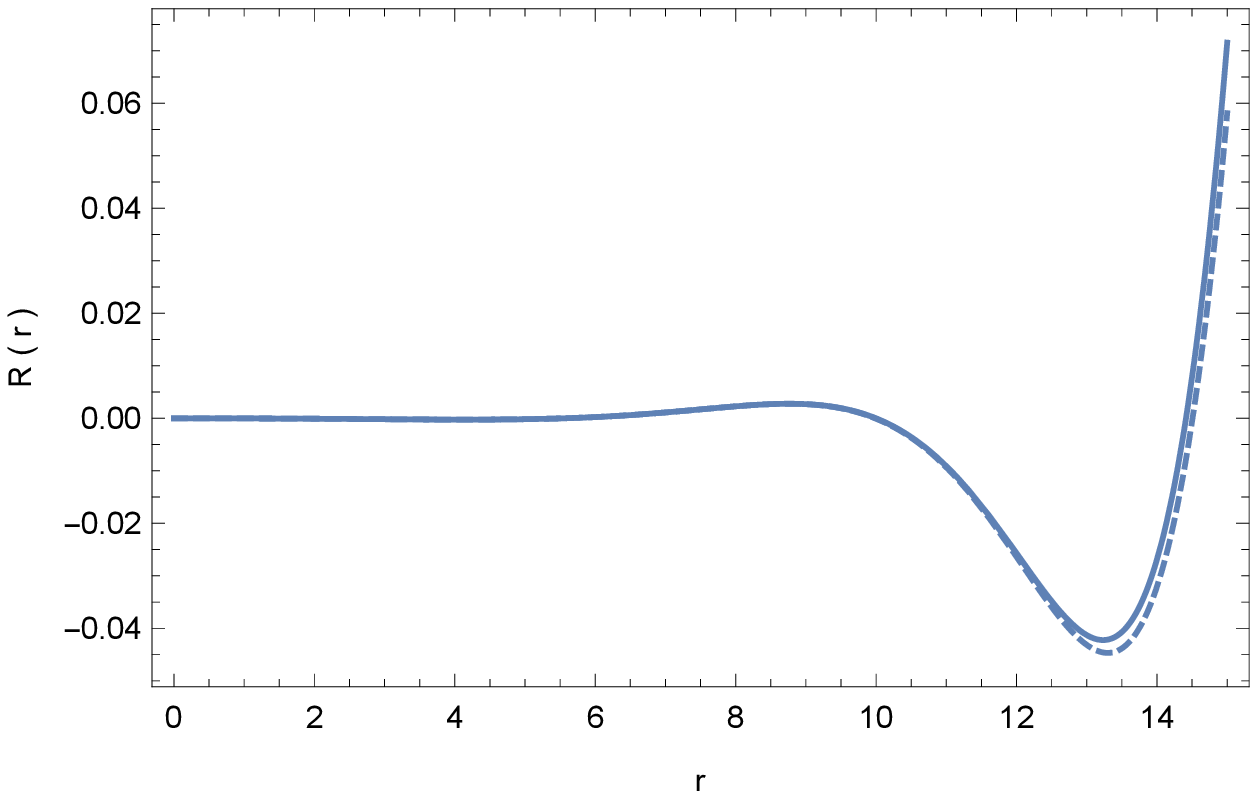} %
\includegraphics[width=8.5cm]{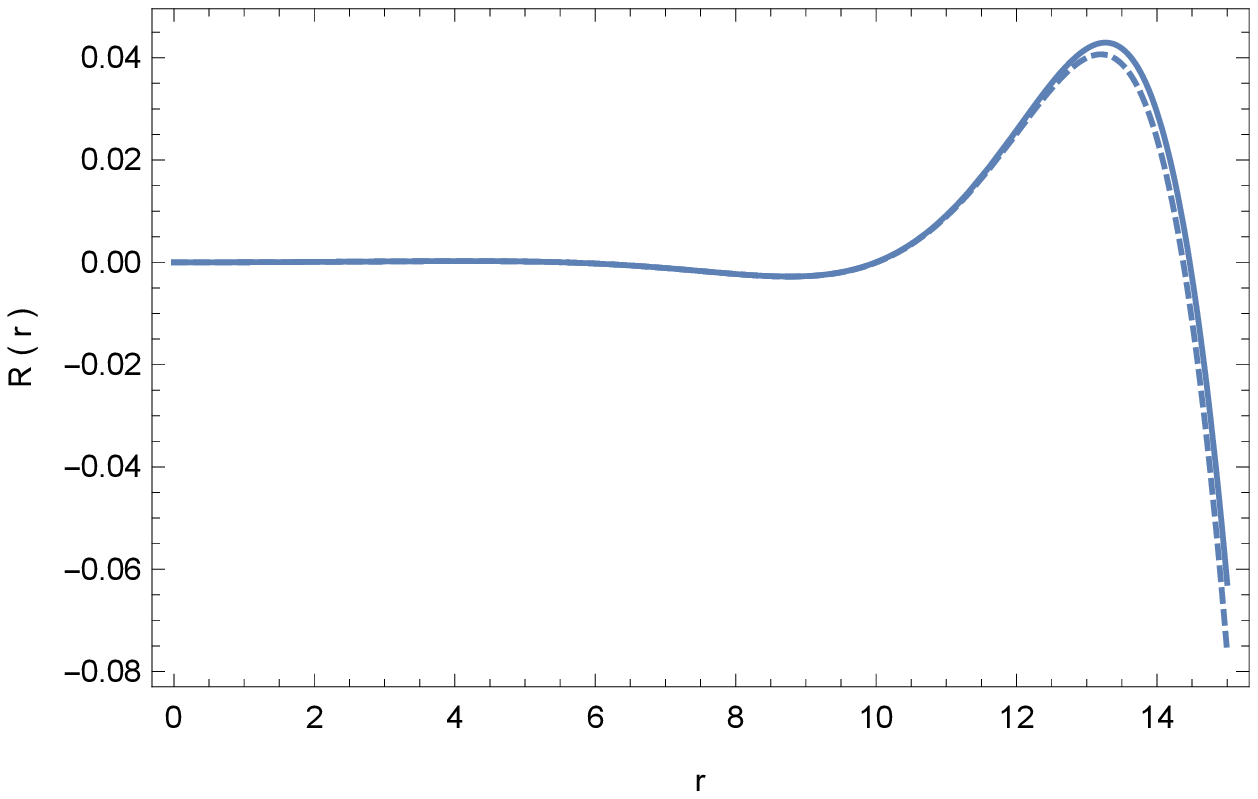}
\caption{Comparison of the numerical solutions of the radial biharmonic
standing wave equation (\protect\ref{sw}) and of the approximate solutions
obtained by the Laplace-Adomian Decomposition Method, truncated to two
terms, $R(r)\approx R_0(r0+R_1(r)$. The numerical solutions are represented
by the solid curves, while the dashed curves depicts the Laplace-Adomian
Decomposition Method two terms solutions. The initial conditions used to
integrate the equations are $R(0)=-7.85\times 10^{-12}$ and $R^{\prime
\prime}(0)=-4.31\times 10^{-5}$ (left figure) and $R(0)=1.27\times 10^{-12}$
and $R^{\prime \prime}(0)=4.31\times 10^{ -5}$ (right figure), respectively.}
\label{fig3}
\end{figure*}

The comparison of the two terms truncated Laplace-Adomian Decomposition
Method solution, $R(r)\approx R_{0}(r)+R_1(r)$ with the exact numerical
solution is presented, for two different sets of initial conditions, in Fig.~%
\ref{fig3}.

\subsection{The Adomian Decomposition Method for the radial biharmonic
standing waves equation}

We consider now the use of the Adomian Decomposition Method for obtaining a
semi-analytical solution of the radial biharmonic standing waves equation.
For the sake of generality we will consider a more general equation of the
form
\begin{equation}
\frac{d^{4}R}{dr^{4}}+f\left( r\right) \frac{d^{3}R}{dr^{3}}=R^{2\sigma
+1}-R,  \label{95}
\end{equation}%
where $f(r)$ is an arbitrary function of the radial coordinate $r$, and
which we will solve with the initial conditions $R(0)\neq 0$, $R^{\prime
}\left( 0\right) =0$, $R^{\prime \prime }\left( 0\right) \neq 0$, and $%
R^{\prime \prime \prime }(0)=0$, respectively. Then the following identity
can be immediately obtained,
\begin{eqnarray}
&&\int_{0}^{r}dr_{1}\int_{0}^{r_{1}}dr_{2}\int_{0}^{r_{2}}e^{-\int f\left(
r_{3}\right) dr_{3}}dr_{3}\int_{0}^{r_{3}}e^{\int f\left( r_{4}\right)
dr_{4}}\left[ R^{\prime \prime \prime \prime }\left( r_{4}\right) +f\left(
r_{4}\right) R^{\prime \prime \prime }\left( r_{4}\right) \right] dr_{4}
\notag \\
&=&\int_{0}^{r}dr_{1}\int_{0}^{r_{1}}dr_{2}\int_{0}^{r_{2}}e^{-\int f\left(
r_{3}\right) dr_{3}}dr_{3}\left[ \int_{0}^{r_{3}}e^{\int f\left(
r_{4}\right) dr_{4}}dR^{\prime \prime \prime \prime}\left( r_{4}\right)
+\int_{0}^{r_{3}}e^{\int f\left( r_{4}\right) dr_{4}}f\left( r_{4}\right)
R^{\prime \prime \prime }\left( r_{4}\right) dr_{4}\right]  \notag \\
&=&\int_{0}^{r}dr_{1}\int_{0}^{r_{1}}dr_{2}\int_{0}^{r_{2}}e^{-\int f\left(
r_{3}\right) dr_{3}}\left[ \int_{0}^{r_{3}}d\left( e^{\int f\left( r\right)
dr}R^{\prime \prime \prime }\left( r\right) \right) \right] dr_{3}  \notag \\
&=&\int_{0}^{r}dr_{1}\int_{0}^{r_{1}}dr_{2}\int_{0}^{r_{2}}e^{-\int f\left(
r_{3}\right) dr_{3}}\left.\left( e^{\int f\left( r\right) dr}R^{\prime
\prime \prime }\left( r\right) \right) \right|_{r=0}^{r=r_{3}}dr_{3}  \notag
\\
&=&\int_{0}^{r}dr_{1}\int_{0}^{r_{1}}dr_{2}\int_{0}^{r_{2}}e^{-\int f\left(
r_{3}\right) dr_{3}}e^{\int f\left( r_{3}\right) dr_{3}}R^{\prime \prime
\prime }\left( r_{3}\right)
dr_{3}=\int_{0}^{r}dr_{1}\int_{0}^{r_{1}}dr_{2}\int_{0}^{r_{2}}R^{\prime
\prime \prime }\left( r_{3}\right) dr_{3}  \notag \\
&=&R\left( r\right) -R\left( 0\right) -R^{\prime \prime }\left( 0\right)
\frac{r^{2}}{2}.
\end{eqnarray}

Thus Eq.~(\ref{95}) can be reformulated as an equivalent integral equation
given by
\begin{equation}
R\left( r\right) =R\left( 0\right) +R^{\prime \prime }\left( 0\right) \frac{%
r^{2}}{2}+\int_{0}^{r}dr_{1}\int_{0}^{r_{1}}dr_{2}\int_{0}^{r_{2}}e^{-\int
f\left( r_{3}\right) dr_{3}}dr_{3}\int_{0}^{r_{3}}e^{\int f\left(
r_{4}\right) dr_{4}}\left[ R^{2\sigma +1}\left( r_{4}\right) -R\left(
r_{4}\right) \right] dr_{4}.
\end{equation}%
By taking into account that $f(r)=4/r$, and by decomposing $R$ and $%
R^{2\sigma +1}$ as $R=\sum_{n=0}^{\infty }{R_{n}}$ and $R^{2\sigma
+1}=\sum_{n=0}^{\infty }{A_{n}}$, where $A_{n}$ are the Adomian polynomials,
we obtain
\begin{equation}
R_{0}\left( r\right) +\sum_{n=0}^{\infty }R_{n+1}\left( r\right) =R\left(
0\right) +R^{\prime \prime }\left( 0\right) \frac{r^{2}}{2}%
+\int_{0}^{r}dr_{1}\int_{0}^{r_{1}}dr_{2}\int_{0}^{r_{2}}\frac{1}{r_{3}^{4}}%
dr_{3}\int_{0}^{r_{3}}r_{4}^{4}\left[ \sum_{n=0}^{\infty }A_{n}\left(
r_{4}\right) -\sum_{n=0}^{\infty }R_{n}\left( r_{4}\right) \right] dr_{4}.
\end{equation}

Then an analytic solution to Eq.~(\ref{95}) can be obtained with the help of
the recursive relations
\begin{equation}
R_{0}\left( r\right) =R\left( 0\right) +R^{\prime \prime }\left( 0\right)
\frac{r^{2}}{2},
\end{equation}%
\begin{equation}
R_{k+1}\left( r\right)
=\int_{0}^{r}dr_{1}\int_{0}^{r_{1}}dr_{2}\int_{0}^{r_{2}}\frac{1}{r_{3}^{4}}%
dr_{3}\int_{0}^{r_{3}}r_{4}^{4}\left[ A_{k}\left( r_{4}\right) -R_{k}\left(
r_{4}\right) \right] dr_{4}.
\end{equation}

\subsubsection{Application: the case $\protect\sigma =1/2$}

As an application of the Adomian Decomposition Method for obtaining the
solution of Eq.~(\ref{95}) we consider the case $\sigma =1/2$. Hence the
radial biharmonic standing wave equation becomes
\begin{equation}  \label{101}
\frac{d^{4}R}{dr^{4}}+\frac{4}{r} \frac{d^{3}R}{dr^{3}}=R^{2}-R.
\end{equation}

The first few terms in the series solution of this equation are given by
\begin{equation}
R_{0}\left( r\right) =R\left( 0\right) +R^{\prime \prime }\left( 0\right)
\frac{r^{2}}{2},
\end{equation}
\begin{equation}
R_1(r)=\frac{1}{120} \left[R(0)-1\right]R(0)r^4+\frac{\left[2 R(0)-1\right]
R^{\prime \prime }(0)}{1680}r^6+\frac{\left(R^{\prime \prime }(0)\right)^2}{%
12096}r^8,
\end{equation}
\begin{eqnarray}
R_2(r)&=&\frac{R(0)\left[2 R(0)^2-3 R(0)+1\right]}{362880}r^8+\frac{ \left[%
18 R(0)^2-18 R(0)+1\right] R^{\prime \prime }(0)}{13305600}r^{10}+\frac{41%
\left[2 R(0)-1\right]\left(R^{\prime \prime }(0)\right)^2}{1037836800}r^{12}+
\notag \\
&& \frac{\left(R^{\prime \prime }(0)\right)^3}{396264960}r^{14},
\end{eqnarray}
\begin{eqnarray}
R_3(r)&=&\frac{R(0) \left[146 R(0)^3-292 R(0)^2+151 R(0)-5\right]}{%
31135104000}r^{12}+ \frac{ \left[1120 R(0)^3-1680 R(0)^2+566 R(0)-3\right]%
R^{\prime \prime }(0)}{1307674368000}r^{14}+  \notag \\
&& \frac{ \left[31282 R(0)^2-31282 R(0)+3407\right] \left(R^{\prime \prime
}(0)\right)^2}{414968666112000}r^{16}+\frac{3061 \left[2 R(0)-1\right]
\left(R^{\prime \prime }(0)\right)^3}{2027418340147200}r^{18}+  \notag \\
&&\frac{89\left(R^{\prime \prime }(0)\right)^4}{1366067972505600} r^{20}.
\end{eqnarray}

\begin{figure*}[htp!]
\centering
\includegraphics[width=8.5cm]{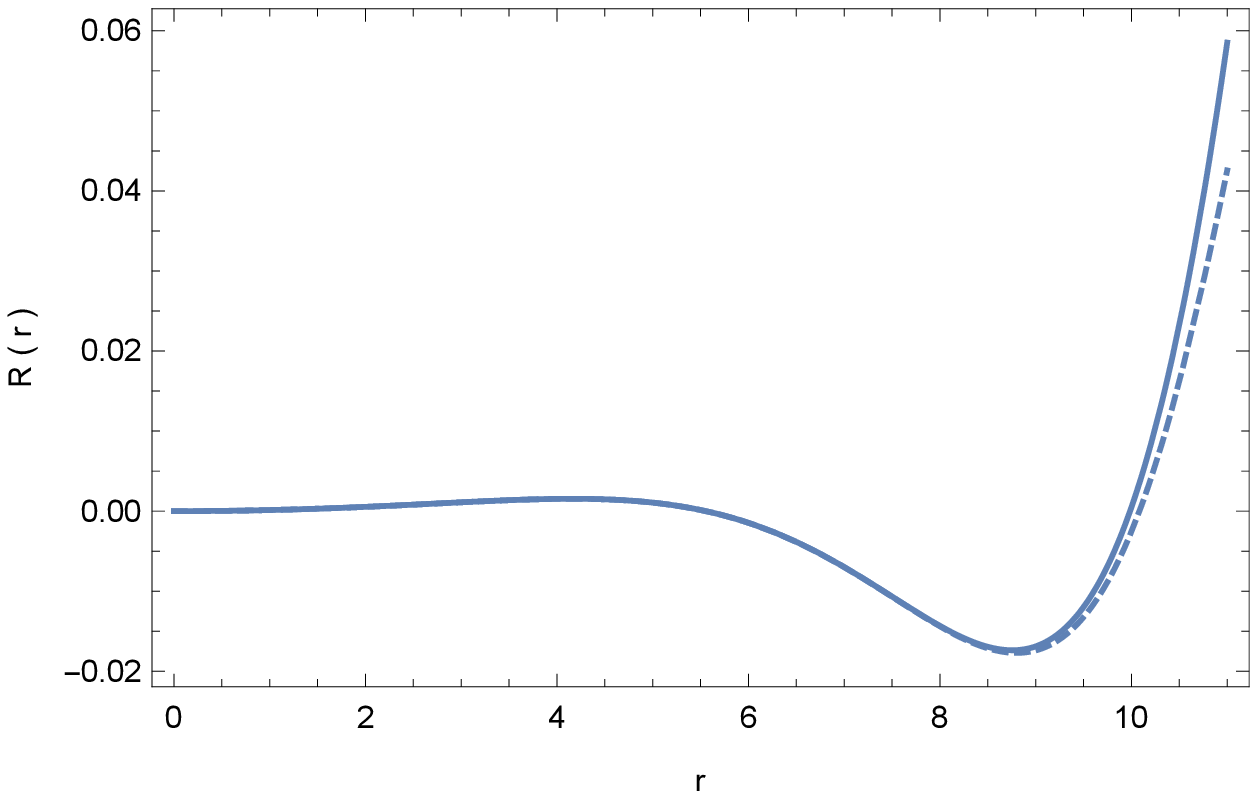} %
\includegraphics[width=8.5cm]{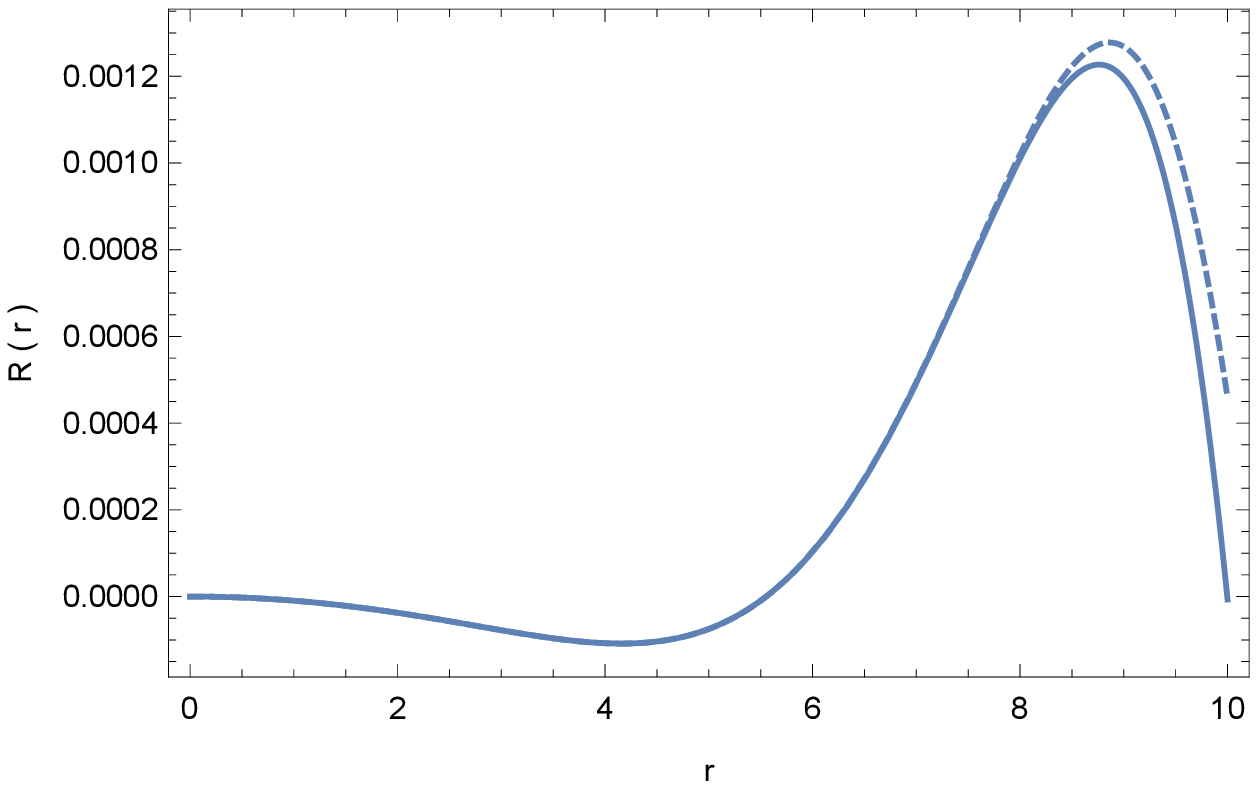}
\caption{Comparison of the numerical solutions of the radial biharmonic
standing wave equation (\protect\ref{101}) and of the Adomian Decomposition
Method approximate solutions, truncated to four terms. The numerical
solutions are represented by the solid curves, while the dashed curves
depicts the Adomian Decomposition Method four terms solutions. The initial
conditions used to integrate the equations are $R(0)=7.44\times 10^{-15}$
and $R^{\prime \prime}(0)=2.71\times 10^{-4}$ (left figure) and $%
R(0)=-3.89\times 10^{-16}$ and $R^{\prime \prime}(0)=-1.91\times 10^{ -5}$
(right figure), respectively.}
\label{fig4}
\end{figure*}

The comparison between the exact numerical solution and the approximate
solution $R(r)=R_0(r)+R_1(r)+R_2(r)+R_3(r)$ of Eq.~(\ref{101}) is
represented, for two sets of initial values, in Fig.~\ref{fig4}.

\section{Discussions and concluding remarks}

\label{sect4}

In the present paper we have presented the applications of the Adomian
Decomposition method for solving the nonlinear biharmonic differential
equation. The Adomian Decomposition Method has been successfully used to
solve many classes of differential, integral and functional equations. It
has also important applications in science and engineering. The basic
ingredient of this approach is the decomposition of the nonlinear term in
the differential equations into a series of polynomials of the form $\sum
_{n=1}^{\infty}{A_n}$, where $A_n$ are the so-called Adomian polynomials.
Simple formulas that can generate Adomian polynomials for many forms of
nonlinearity have been derived in \cite{R1,R2}. The solutions of the
nonlinear differential equations can be obtained recursively, and each term
of the Adomian series can be computed once the corresponding polynomial,
obtained from an expansion of the nonlinear term into a power series, is
known.

We have considered in detail both the one dimensional, as well as the
radial, three dimensional, biharmonic type equation containing some
nonlinear terms. We have implemented two versions of the Adomian
Decomposition Method for solving the biharmonic equation, namely, the
Laplace-Adomian Decomposition Method, and the standard Adomian Decomposition
Method. The Laplace-Adomian Decomposition Method combines the powerful
Laplace transformation with the advantages of the Adomian method, with the
iterative procedure applied in the space of the Laplace transformed
functions. In the radial case the Laplace transforms of the terms in the
Adomian expansion can be obtained as solutions of a first order differential
equation, which can be obtained by quadratures. However, in the present case
the integral, and the Laplace transform itself, cannot be obtained in an
exact form, and therefore one have to resort to some approximate methods.

For each type of considered equations we have also considered some concrete
examples, and we have compared the Adomian solution with the exact numerical
solution. Generally, the efficiency and the precision of the Laplace-Adomian
Decomposition Method is very good. In the case of the one-dimensional
standing wave biharmonic equation only three terms of the Adomian expansion
are enough to give a good approximation of the numerical solution, while for
the case of the radial nonlinear biharmonic standing wave equation the
numerical solution can be approximated by using only two terms. This
coincidence implicitly shows the power of the Adomian method, which can be
used to find out even the exact solution of a given differential equation.
However, in general the application of the method may be complicated by the
difficulties in solving exactly the differential equations for the Laplace
transform, and for obtaining the inverse Laplace transform. But, at least in
the case of the radial nonlinear standing wave equation, a simple technique
based on the power series expansion of the Laplace transform of the Adomian
polynomials gives good approximations of the numerical solutions. Numerical
techniques for obtaining the inverse Laplace transform \cite{Lern} may also
be useful in obtaining the successive terms in the Laplace-Adomian expansion.

We have also considered the standard Adomian Decomposition Method for both
the first order and radial nonlinear biharmonic equations. Computationally,
this method is very simple, and it can provide some power series solutions
that can describe the numerical solution relatively well. The Adomian method
is very simple and efficient, but it may raise some questions about the
convergence of the series of functions \cite{C1,C2}. Moreover, we must point
out that the accuracy of the approximations of the numerical solution by the
Adomian series is strongly dependent on the initial conditions used to solve
the equations. The Adomian solutions work well for small numerical values of
the initial conditions. Once these values are increased, the accuracy of the
estimations becomes poor, at least for the number of terms used to
approximate the solutions in the present approach. These raises the issue of
the dependence of the convergence of the Adomian solution from the initial
conditions, a mathematical problem certainly worth of investigating.

The biharmonic equation appears in many physical and engineering applications \cite{Bose1,Boos, n1, n2}. In particular, it plays an important role within the hydrodynamic formulation of the Schr\"{o}dinger equation, and in the presence of the quantum potential. This physical approach is extensively used for the study of the quantum fluids. In many applications, mostly due to the computational difficulties, the quantum potential is neglected. However, by using the present approach, semi-analytical solutions of the biharmonic equation can be obtained, which can approximate well the numerical solution. The semi-analytical solutions offer the possibility of a deeper insight into the physical nature of the problem, as well as of a significant simplification of the estimation of some relevant physical parameters. Similar applications of the method could lead to the development of powerful mathematical methods for solving different problems described by fourth order differential equations that play an important role in engineering, like, for example, in the study of the large amplitude free vibrations of a uniform cantilever beam \cite{Bel}.

The Adomian Decomposition Method, as well as its Laplace transform version,
represents a powerful mathematical tool for physicists and engineers
investigating both theoretical and applied problems. The biharmonic
equation, and its extensions, are interesting in themselves from a
mathematical point of view. There are also important in many applications.
In the present study we have introduced some theoretical tools, which are
extremely effective in dealing with strongly nonlinear differential
equations and complex mathematical models, and that may help in the better
understanding of the properties and solutions of the nonlinear biharmonic
equation.

\section*{Acknowledgments}

We would like to thank the anonymous reviewer for comments and suggestions that helped us to improve our manuscript. T. H. would like to thank the Yat Sen School of the Sun Yat Sen University in Guangzhou for the kind hospitality offered during the preparation of this work.

\end{document}